\documentclass[12pt,twoside,leqno]{article}
\usepackage{amsmath}
\usepackage{amssymb}
\usepackage{amsxtra}
\usepackage{amscd}
\usepackage{amsthm}
\usepackage[mathscr]{eucal}

\theoremstyle{plain}
\newtheorem{thm}[subsection]{Theorem}

\newtheorem{lem}[subsection]{Lemma}

\newtheorem{sbprop}[subsubsection]{Proposition}

\newtheorem{sblem}[subsubsection]{Lemma}

\theoremstyle{definition}

\newtheorem{para}[subsection]{}

\newtheorem{sbeg}[subsubsection]{Example}

\begin{document}

\title{A description of a result of Deligne\\ 
by log higher Albanese map}

\author
{{\normalsize Dedicated to Goo Ishikawa on his sixtieth birthday.}\\
\\
Sampei Usui}

\maketitle
\renewcommand{\mathbb}{\bold}

\newcommand\Cal{\mathcal}
\newcommand\define{\newcommand}
\define\gp{\mathrm{gp}}%
\define\fs{\mathrm{fs}}%
\define\an{\mathrm{an}}%
\define\mult{\mathrm{mult}}%
\define\add{\mathrm{add}}%
\define\Ker{\mathrm{Ker}\,}%
\define\Coker{\mathrm{Coker}\,}%
\define\Hom{\mathrm{Hom}\,}%
\define\Ext{\mathrm{Ext}\,}%
\define\rank{\mathrm{rank}\,}%
\define\gr{\mathrm{gr}}%
\define\cHom{\Cal{Hom}}
\define\cExt{\Cal Ext\,}%

\define\cC{\Cal C}
\define\cD{\Cal D}
\define\cO{\Cal O}
\define\cS{\Cal S}
\define\cM{\Cal M}
\define\cG{\Cal G}
\define\cH{\Cal H}
\define\cE{\Cal E}
\define\cF{\Cal F}
\define\cN{\Cal N}
\define\cQ{\Cal Q}
\define\fF{\frak F}
\define\fg{\frak g}
\define\fh{\frak h}
\define\Dc{\check{D}}
\define\Ec{\check{E}}

\newcommand{\N}{{\mathbb{N}}}
\newcommand{\Q}{{\mathbb{Q}}}
\newcommand{\Z}{{\mathbb{Z}}}
\newcommand{\R}{{\mathbb{R}}}
\newcommand{\C}{{\mathbb{C}}}
\newcommand{\bN}{{\mathbb{N}}}
\newcommand{\bQ}{{\mathbb{Q}}}
\newcommand{\bF}{{\mathbb{F}}}
\newcommand{\bZ}{{\mathbb{Z}}}
\newcommand{\bP}{{\mathbb{P}}}
\newcommand{\bR}{{\mathbb{R}}}
\newcommand{\bC}{{\mathbb{C}}}
\newcommand{\bS}{{\bold{S}}}
\newcommand{\bbQ}{{\bar \mathbb{Q}}}
\newcommand{\ol}[1]{\overline{#1}}
\newcommand{\too}{\longrightarrow}
\newcommand{\respect}{\rightsquigarrow}
\newcommand{\compatible}{\leftrightsquigarrow}
\newcommand{\upc}[1]{\overset {\lower 0.3ex \hbox{${\;}_{\circ}$}}{#1}}
\newcommand{\Gmlog}{\bG_{m, \log}}
\newcommand{\Gm}{\bG_m}
\newcommand{\ep}{\varepsilon}
\newcommand{\Spec}{\operatorname{Spec}}
\newcommand{\val}{{\mathrm{val}}} 
\newcommand{\n}{\operatorname{naive}}
\newcommand{\bs}{\operatorname{\backslash}}
\newcommand{\Gal}{\operatorname{{Gal}}}
\newcommand{\gal}{{\rm {Gal}}({\bar \Q}/{\Q})}
\newcommand{\galp}{{\rm {Gal}}({\bar \Q}_p/{\Q}_p)}
\newcommand{\gall}{{\rm{Gal}}({\bar \Q}_\ell/\Q_\ell)}
\newcommand{\wep}{W({\bar \Q}_p/\Q_p)}
\newcommand{\wel}{W({\bar \Q}_\ell/\Q_\ell)}
\newcommand{\Ad}{{\rm{Ad}}}
\newcommand{\BS}{{\rm {BS}}}
\newcommand{\even}{\operatorname{even}}
\newcommand{\End}{{\rm {End}}}
\newcommand{\odd}{\operatorname{odd}}
\newcommand{\GL}{\operatorname{GL}}
\newcommand{\np}{\text{non-$p$}}
\newcommand{\g}{{\gamma}}
\newcommand{\G}{{\Gamma}}
\newcommand{\Lam}{{\Lambda}}
\newcommand{\La}{{\Lambda}}
\newcommand{\lam}{{\lambda}}
\newcommand{\la}{{\lambda}}
\newcommand{\uL}{{{\hat {L}}^{\rm {ur}}}}
\newcommand{\uQp}{{{\hat \Q}_p}^{\text{ur}}}
\newcommand{\sel}{\operatorname{Sel}}
\newcommand{\dt}{{\rm{Det}}}
\newcommand{\Sig}{\Sigma}
\newcommand{\fil}{{\rm{fil}}}
\newcommand{\SL}{{\rm{SL}}}
\newcommand{\spl}{{\rm{spl}}}
\newcommand{\st}{{\rm{st}}}
\newcommand{\Isom}{{\rm {Isom}}}
\newcommand{\Mor}{{\rm {Mor}}}
\newcommand{\bg}{\bar{g}}
\newcommand{\id}{{\rm {id}}}
\newcommand{\cone}{{\rm {cone}}}
\newcommand{\al}{a}
\newcommand{\ChL}{{\cal{C}}(\La)}
\newcommand{\Image}{{\rm {Image}}}
\newcommand{\toric}{{\operatorname{toric}}}
\newcommand{\torus}{{\operatorname{torus}}}
\newcommand{\Aut}{{\rm {Aut}}}
\newcommand{\Qp}{{\mathbb{Q}}_p}
\newcommand{\barQp}{{\mathbb{Q}}_p}
\newcommand{\Qpur}{{\mathbb{Q}}_p^{\rm {ur}}}
\newcommand{\Zp}{{\mathbb{Z}}_p}
\newcommand{\Zl}{{\mathbb{Z}}_l}
\newcommand{\Ql}{{\mathbb{Q}}_l}
\newcommand{\Qlur}{{\mathbb{Q}}_l^{\rm {ur}}}
\newcommand{\F}{{\mathbb{F}}}
\newcommand{\eps}{{\epsilon}}
\newcommand{\epsLa}{{\epsilon}_{\La}}
\newcommand{\epsLaVxi}{{\epsilon}_{\La}(V, \xi)}
\newcommand{\epsOLaVxi}{{\epsilon}_{0,\La}(V, \xi)}
\newcommand{\Qplin}{{\mathbb{Q}}_p(\mu_{l^{\infty}})}
\newcommand{\otimesQplin}{\otimes_{\Qp}{\mathbb{Q}}_p(\mu_{l^{\infty}})}
\newcommand{\galFl}{{\rm{Gal}}({\bar {\Bbb F}}_\ell/{\Bbb F}_\ell)}
\newcommand{\gallur}{{\rm{Gal}}({\bar \Q}_\ell/\Q_\ell^{\rm {ur}})}
\newcommand{\galFF}{{\rm {Gal}}(F_{\infty}/F)}
\newcommand{\galFv}{{\rm {Gal}}(\bar{F}_v/F_v)}
\newcommand{\galF}{{\rm {Gal}}(\bar{F}/F)}
\newcommand{\epsVxi}{{\epsilon}(V, \xi)}
\newcommand{\epsOVxi}{{\epsilon}_0(V, \xi)}
\newcommand{\plim}{\lim_
{\scriptstyle 
\longleftarrow \atop \scriptstyle n}}
\newcommand{\sig}{{\sigma}}
\newcommand{\ga}{{\gamma}}
\newcommand{\del}{{\delta}}
\newcommand{\Vss}{V^{\rm {ss}}}
\newcommand{\Bst}{B_{\rm {st}}}
\newcommand{\Dpst}{D_{\rm {pst}}}
\newcommand{\Dcrys}{D_{\rm {crys}}}
\newcommand{\DdR}{D_{\rm {dR}}}
\newcommand{\Fin}{F_{\infty}}
\newcommand{\Kla}{K_{\lambda}}
\newcommand{\Ola}{O_{\lambda}}
\newcommand{\Mla}{M_{\lambda}}
\newcommand{\Det}{{\rm{Det}}}
\newcommand{\Sym}{{\rm{Sym}}}
\newcommand{\LaSa}{{\La_{S^*}}}
\newcommand{\cX}{{\cal {X}}}
\newcommand{\MHG}{{\frak {M}}_H(G)}
\newcommand{\tauMla}{\tau(M_{\lambda})}
\newcommand{\Fvur}{{F_v^{\rm {ur}}}}
\newcommand{\Lie}{{\rm {Lie}}\,}
\newcommand{\cB}{{\cal {B}}}
\newcommand{\cL}{{\cal {L}}}
\newcommand{\cW}{{\cal {W}}}
\newcommand{\fq}{{\frak {q}}}
\newcommand{\cont}{{\rm {cont}}}
\newcommand{\SC}{{SC}}
\newcommand{\Om}{{\Omega}}
\newcommand{\dR}{{\rm {dR}}}
\newcommand{\crys}{{\rm {crys}}}
\newcommand{\hatSig}{{\hat{\Sigma}}}
\newcommand{\rdet}{{{\rm {det}}}}
\newcommand{\ord}{{{\rm {ord}}}}
\newcommand{\Alb}{{{\rm {Alb}}}}
\newcommand{\BdR}{{B_{\rm {dR}}}}
\newcommand{\BdRO}{{B^0_{\rm {dR}}}}
\newcommand{\Bcrys}{{B_{\rm {crys}}}}
\newcommand{\Qw}{{\mathbb{Q}}_w}
\newcommand{\barkappa}{{\bar{\kappa}}}
\newcommand{\cP}{{\Cal {P}}}
\newcommand{\cZ}{{\Cal {Z}}}
\newcommand{\oppLa}{{\Lambda^{\circ}}}
\newcommand{\bG}{{\mathbb{G}}}
\newcommand{\br}{{{\bold r}}}
\newcommand{\triv}{{\rm{triv}}}
\newcommand{\sub}{{\subset}}
\newcommand{\LD}{{D^{\star,\mild}_{\SL(2)}}}
\newcommand{\LbD}{{D^{\star}_{\SL(2)}}}
\newcommand{\dbDv}{{D^{\star}_{\SL(2),\val}}}
\newcommand{\nspl}{{{\rm nspl}}}
\newcommand{\lval}{{[\val]}}
\newcommand{\mild}{{{\rm{mild}}}}
\newcommand{\lan}{\langle}
\newcommand{\ran}{\rangle}
\newcommand{\rar}{{{\rm {rar}}}}
\newcommand{\rac}{{{\rm {rac}}}}
\newcommand{\Rep}{{\operatorname{Rep}}}
\newcommand{\LMH}{{\rm {LMH}}}
\newcommand{\MHS}{{\rm {MHS}}}
\newcommand{\red}{{{\rm red}}}
\newcommand{\cV}{{\Cal {V}}}

\begin{abstract}

\noindent
In a joint work \cite{KNU17} with Kazuya Kato and Chikara Nakayama, 
log higher Albanese manifolds were constructed as an application of log mixed Hodge theory with group action.
In this framework, we describe a work of Deligne in \cite{D89} on some nilpotent quotients of the fundamental group of the projective line minus three points, where polylogarithms appear.
As a result, we have $q$-expansions of higher Albanese maps at boundary points, i.e., log higher Albanese maps over the boundary.

\end{abstract}


\section*{Contents}

\noindent \S\ref{s:intro}. Introduction

\noindent \S\ref{s:glmh}. Log mixed Hodge structures with group action

\noindent \S\ref{s:lha}. Log higher Albanese manifolds

\noindent \S\ref{s:example}. Description of a result of Deligne by log higher Albanese map

\noindent \S\ref{s:appen}. Summary of a result of Deligne in \cite{D89}

\renewcommand{\thefootnote}{\fnsymbol{footnote}}
\footnote[0]{Primary 14C30; 
Secondary 14D07, 32G20} 
\renewcommand{\thefootnote}{\arabic{footnote}}

\setcounter{section}{-1}
\section{Introduction}\label{s:intro}

We review the results of \cite{KNU17}:
- General theory of log mixed Hodge structures with polarizable graded quotients endowed with group actions.
- Description of the functors represented by higher Albanese manifolds in terms of tensor functors.
- Toroidal partial compactifications of higher Albanese manifolds to get log higher Albanese manifolds, and describe the functors represented by them.

We describe a result of Deligne in \cite{D89} about polylogarithms, which were appeared in higher Albanese maps, in terms of the log higher Albanese maps.
The advantage of our formulation is that log higher Albanese maps have $q$ expansions at the boundary points over which we observe directly $\zeta(n)$ $(n\ge2)$ as values of polylogarithms.

For readers' convenience, we add as Appendix a summary of the related result of Deligne in \cite{D89}.

Actually, for the present description in Section 3, it is enough to use the formulation of spaces of nilpotent orbits in \cite{KNU} Part III.
The formulation of them in \cite{KNU17} is reviewed in Sections 1 and 2 for further study in the case of higher Albanese manifolds with non-trivial Hodge structures.

\section{Log mixed Hodge structures with group action}\label{s:glmh}

\medskip

We review general formulations and results of log mixed Hodge structures with group action in \cite{KNU17} and \cite{KNU} Parts III, IV, in a  minimal size for the later use of this paper.
The full version will be appeared in \cite{KNU} Part V.

\begin{para}\label{log}

A {\it log structure} on a ringed space $(S,\cO_S)$ consists of a sheaf of monoids $M$ on $S$ and a homomorphism $\alpha:M\to\cO_S$ such that $\alpha^{-1}(\cO_S^\times)\overset\sim\to\cO_S^\times$.
\medskip

\begin{sbeg}
Let $S=\bC$ and $\{0\}$ a divisor. 
The associated log structure is 
$M=\{ f\in\cO_S\,|\,\text{$f$ is invertible on $S\smallsetminus \{0\}$}\}$.
$S^{\log}$ is defined to be the set of all pairs $(s,h)$ consisting of a point $s\in S$ and an argument function $h$ which is a homomorphism $M_s\to \bS^1$ of monoids whose restriction to $\cO_{S,s}^{\times}$ is $u\mapsto u(s)/|u(s)|$.

In this case, the ringed space $(S^{\log},\cO_S^{\log})$ is explained as follows.
Let $\tilde S^{\log}:=\bC\cup(\bR\times i\infty)=\bR\times i(\bR\cup\infty)$ endowed with coordinate function $z=x+iy$ $(-\infty< y\le\infty)$.
Let $S^{\log}:=(\bC\cup(\bR\times i\infty))/\bZ$, and consider maps
$\tilde S^{\log}\to S^{\log}\to S:$ $z=x+iy\mapsto(e^{-2\pi y}, e^{2\pi ix})\mapsto q:=e^{2\pi iz}$.
Note that $(e^{-2\pi y}, e^{2\pi ix})$ is a polar coordinate extended over $-\infty< y\le\infty$, and $S^{\log}\to S$ is a real oriented blowing-up at $\{0\}$, which is proper.
$h:M_0\to\bold S^1$ in $t:=(0,h)\in S^{\log}$ sends $q$ to $e^{2\pi ix}$.
Since $z$ is considered as a branch of $(2\pi i)^{-1}\log(q)$, we have $\cO_{S,t}^{\log}=\cO_{S,0}[z]$ which is isomorphic to a polynomial algebra $\cO_{S,0}[T]$ of one indeterminate $T$ over $\cO_{S,0}$ under $z\leftrightarrow T$ (\cite{KU09} 2.2.5).

\end{sbeg}

For more general and finer treatment, see \cite{KN99}, \cite{KU09} 2.2.

\end{para}

\begin{para}\label{Gsetting}
Let  $G$ be a linear algebraic group over $\bQ$.
Let $G_u$ be the unipotent radical of $G$ and let $G_{\red}=G/G_u$.
Let $\Rep(G)$ be the category of finite-dimensional linear representations of $G$ over $\bQ$.

\end{para}

\begin{para}\label{k0}
Let $k_0:\bG_m\to G_{\red}$ be a $\Q$-rational and central homomorphism.
Assume that, for one (hence all) lifting $\bG_{m,\bR}\to G_\bR$ of $k_0$, the adjoint action of $\bG_{m,\bR}$ on $\Lie(G_u)_\bR=\bR\otimes_\bQ\Lie(G_u)$ is of weight $\le-1$.

Then, for any $V\in \Rep(G)$, the action of $\bG_m$ on $V$ via a lifting $\bG_m\to G$ of $k_0$ defines an increasing filtration $W$ on $V$ over $\bQ$, called {\it weight filtration}, which is independent of the lifting.
\end{para}

\begin{para}\label{type k0}

Assume that we are given a homomorphism  $k_0:\bG_m\to G_{\red}$ as in \ref{k0}.
A {\it $G$-mixed Hodge structure} ({\it $G$-MHS}, for short) {\it of type $k_0$} is an exact $\otimes$-functor (\cite{D90} 2.7) from $\Rep(G)$ to the category of $\bQ$-mixed Hodge structures keeping the underlying vector spaces  with weight filtrations.

\end{para}

\begin{para}\label{SCR}

As in \cite{D71}, let $S_{\bC/\bR}$ be the Weil restriction of scalars of $\bG_m$ from $\bC$ to $\bR$.
It represents the functor $A\mapsto(\bC\otimes_\bR A)^\times$ for commutative rings $A$ over $\bR$.
In particular, $S_{\bC/\bR}(\bR)=\bC^\times$, which is understood as $\bC^\times$ regarded as an algebraic group over $\bR$.

Let $w:\bG_m\to S_{\bC/\bR}$ be the homomorphism induced from the natural map
$A^\times\to(\bC\otimes_\bR A)^\times$.

\end{para}

\begin{para}\label{hd}

The following (1) and (2) are equivalent:

(1) A finite-dimensional linear representation of $S_{\bC/\bR}$ over $\bR$.

(2) A finite-dimensional $\bR$-vector space $V$ with a decomposition 
$V_\bC:=\bC\otimes_\bR V=\bigoplus_{p,q\in\bZ}V_\bC^{p,q}$ such that, for any $p,q$, $V_\bC^{q,p}$ is complex conjugate of $V_\bC^{p,q}$ (Hodge decomposition).

For a finite-dimensional linear representation $V$ of $S_{\bC/\bR}$, the corresponding decomposition is defined by
$$
V_{\bC}^{p,q}:=\{v\in V_{\bC}\;|\;[z]v=z^p\bar z^qv\;\text{for}\;z\in\bC^\times\}.
$$
Here $[z]$ denotes $z\in\bC^\times$ regarded as an element of $S_{\bC/\bR}(\bR)$.

\end{para}

\begin{para}\label{ah}

Let $H$ be a $G$-MHS of type $k_0$ (\ref{type k0}).
By \ref{hd} and Tannaka duality (cf.\ \cite{D90} 1.12 Th\'eor$\grave{\rm e}$me), the Hodge decompositions of $\gr^W$ of $H(V)$ for $V\in\Rep(G)$ give a homomorphism $S_{\bC/\bR}\to(G_{\red})_\bR$ such that the composite $\bG_m\overset{w}\to S_{\bC/\bR}\to(G_{\red})_\bR$
coincides with $k_0$.
We call this $S_{\bC/\bR}\to(G_{\red})_\bR$ the {\it associated homomorphism with $H$}.

\end{para}

\begin{para}\label{D}

Let $k_0:\bG_m\to G_{\red}$ be as in \ref{k0}.
Fix a homomorphism $h_0:S_{\bC/\bR}\to(G_{\red})_\bR$ such that $h_0\circ w=k_0$.

{\it $G$-mixed Hodge structure of type $h_0$} is a $G$-mixed Hodge structure of type $k_0$ (\ref{type k0}) whose associated homomorphism $S_{\bC/\bR}\to(G_{\red})_\bR$ (\ref{ah}) is $G_{\red}(\bR)$-conjugate to $h_0$.

The {\it period domain $D=D(G,h_0)$} associated to $(G,h_0)$ is defined to be the set of isomorphism classes of $G$-mixed Hodge structures of type $h_0$.

\end{para}

\begin{para}\label{usual}

Usual period domains of Griffiths \cite{G} and their generalization for mixed Hodge structures \cite{U} are special cases of the present period domains.

  Let $\Lambda=(H_0, W, (\langle \;,\,\rangle_w)_w, (h^{p,q})_{p,q})$ be the Hodge data as usual as in \cite{KNU} Part III. 
  Let $G$ be the subgroup of $\Aut(H_{0,\Q}, W)$ consisting of elements which induce {\it similitudes for $\langle\;,\,\rangle_w$} for each $w$. That is, 
$
G:= \{g\in \Aut(H_{0,\Q}, W)\;|\;$ for any $w$, there is a $t_w\in {\mathbb G}_m$ such that 
$\langle gx,gy\rangle_w = t_w\langle x, y \rangle_w$ for any $x,y \in \gr^W_w\}.
$
  Let $G_1:=\Aut(H_{0,\Q}, W, (\langle \;,\,\rangle_w)_w) \subset G$.

  Let $D(\Lambda)$ be the period domain of \cite{U}. 
Then $D(\Lambda)$ is identified with an open and closed part 
 of  $D$ 
in this paper as follows. 

Assume that $D(\Lam)$ is not empty and fix an $\br\in D(\Lambda)$. Then the Hodge decomposition of  $\gr^W\br$  induces 
$h_0: S_{\C/\R} \to (G_{\text{red}})_\R$. (We have
$\langle [z]x, [z]y\rangle_w= |z|^{2w}\langle x, y\rangle_w$ for $z\in \C^\times$ (see \ref{SCR} for $[z]$).) 
  Consider the associated period domain $D$ (\ref{D}). Then $D$ is a finite disjoint union of $G_1(\R)G_u(\C)$-orbits which are open and closed in $D$. 
  Let $\cD$ be the $G_1(\R)G_u(\C)$-orbit in $D$ consisting of points whose associated homomorphisms $S_{\C/\R}\to (D_{\text{red}})_\R$  are $(G_1/G_u)(\R)$-conjugate to $h_0$. 
  Then the map $H\mapsto H(H_{0,\Q})$ gives a $G_1(\R)G_u(\C)$-equivariant  isomorphism $\cD\simeq D(\Lambda)$.

\end{para}

\begin{para}\label{Y}

Fix a homomorphism $h_0:S_{\bC/\bR}\to(G_{\red})_\bR$ as in \ref{D}.

Let $\cC$ be the category of triples $(V,W,F)$ consisting of a finite-dimensional $\bQ$-vector space $V$, an increasing filtration $W$ on $V$ (called the weight filtration), and a decreasing filtration $F$ on $V_\bC$ (called the Hodge filtration).

Let $Y$ be the set of all isomorphism classes of exact $\otimes$-functors from $\Rep(G)$ to $\cC$ preserving the underlying vector spaces with weight filtrations.

Then $G(\bC)$ acts on $Y$ by changing the Hodge filtration $F$.
We have $D\subset Y$ and $D$ is a $G(\bR)G_u(\bC)$-orbit in $Y$ (cf.\ \cite{KNU17} Proposition 3.2.5).
We define $\check D:=G(\bC)D$ in $Y$.
Thus
$$
D\subset \check D=G(\bC)D \subset Y.
$$

$\Dc$ is a $G(\bC)$-homogeneous space and $D$ is an open subspace.
Hence $\Dc$ and $D$ are complex analytic manifolds.

\end{para}


\begin{para}\label{pol} 
  Let $h_0: S_{\C/\R}\to (G_{\red})_\R$ be as in \ref{D}. 
Let $C$ be the image of $i\in \C^\times = S_{\C/\R}(\R)$ by $h_0$ in $(G_{\red})(\R)$ (Cartan involution). 
  We say that $h_0$ is {\it $\R$-polarizable} if 
$\{a\in (G_{\red})'(\R)\;|\; Ca=aC\}$ is a maximal compact subgroup of $(G_{\red})'(\R)$.
Here $(G_{\red})'$ denotes the commutator subgroup of $G_{\red}$.

\end{para}

\begin{para}\label{Gamma1}

Let $h_0: S_{\C/\R}\to (G_{\red})_\R$ be $\R$-polarizable (\ref{pol}).

Let $\Gamma$ be a subgroup of $G(\Q)$ for which there is a faithful $V\in \Rep(G)$ and a $\Z$-lattice $L$ in $V$ such that $L$ is stable under the action of $\Gamma$.

Then the following holds (\cite{KNU17} Proposition 3.3.4):

$(1)$ The action of $\Gamma$ on $D$ is proper, and the quotient space $\Gamma \bs D$ is Hausdorff. 

$(2)$ If $\Gamma$ is torsion-free and if $\gamma p=p$ with $\gamma\in \Gamma$ for  some $p\in D$, then $\gamma=1$.

$(3)$ If $\Gamma$ is torsion-free, then the projection $D\to \Gamma \bs D$ is a local homeomorphism.

\end{para}

\begin{para}\label{nilp}
Let $(G,h_0)$ be as above.

A {\it nilpotent cone} is a cone $\sigma$ over $\bR_{\ge0}$ in $\Lie(G)_\bR$ generated by a finite number of mutually commuting elements such that, for any $V\in\Rep(G)$, the image of $\sigma$ under the induced map $\Lie(G)_\bR\to \End_\bR(V)$ consists of nilpotent operators.

For $F\in\check D$ and a nilpotent cone $\sigma$, $(\sigma,\exp(\sigma_\bC)F)$ is a {\it nilpotent orbit} if it satisfies the following conditions:
Take a generators $N_1,\dots,N_n\in\Lie(G)_\bR$ of the cone $\sigma$.

(1)  (admissibility) There is a faithful $V\in\Rep(G)$ such that the relative monodromy weight filtrations $M(N_j,W)$ on $V$ exist for all $1\le j\le n$.

(2) (Griffiths transversality) $N_jF^p\sub F^{p-1}$ for any $1\le j\le n$,  $p\in \bZ$.

(3) (limit mixed Hodge property) $\exp(\sum_{j=1}^n iy_jN_j)F\in D$ if $y_j \in \bR_{>0}$ are sufficiently large.
\medskip

\end{para}

This is well-defined, i.e., independent of choices of generators $N_1,\dots,N_n$.

Note that, for admissibility, the above condition (1) is enough under the assumption of $\bR$-polarizability (cf.\ \cite{Kas86}, \cite{KNU} III Proposition 1.3.4, Remark in 2.2.2, \cite{Kat14} Proposition 2.1.10).

\begin{para}\label{DSig}

  A {\it weak fan $\Sig$ in $\Lie(G)$} is a nonempty set of sharp rational nilpotent cones satisfying the conditions that it is closed under taking faces and that any $\sig, \sig' \in \Sig$ coincide if they have a common interior point and if there is an $F\in \Dc$ such that both $(\sigma,\exp(\sigma_\bC)F)$ and $(\sigma',\exp(\sigma'_\bC)F)$ are nilpotent orbits. 

For a weak fan $\Sig$ in $\Lie(G)$, let $D_\Sigma$ be the set of all nilpotent orbits $(\sigma,\exp(\sigma_\bC)F)$ with $\sigma\in\Sigma$ and $F\in\Dc$.

\end{para}

\begin{para}\label{sc}

Let $\Gamma$ be a subgroup of $G(\Q)$ as in \ref{Gamma1}.

A weak fan $\Sigma$ in \ref{DSig} is said to be {\it strongly compatible with $\Gamma$} if $\Sigma$ is stable under the adjoint action of $\Gamma$ and each cone $\sigma\in\Sigma$ is generated over $\bR_{\ge0}$ by log of $\exp(\sigma)\cap\Gamma$.

\end{para}

\begin{para}\label{cB}

$\cB$ denotes the category of locally ringed spaces with a covering by open sets each of which has the strong topology in an analytic space.
$\cB(\log)$ denotes the category of objects of $\cB$ endowed with an fs log structure.
For precise definitions of these, see \cite{KU09} 3.2.4, \cite{KNU} Part III 1.1.

\end{para}

\begin{para}\label{LMH}

Let $S\in\cB(\log)$.
A $\bQ$-{\it log mixed Hodge structure} ({\it $\bQ$-LMH}, for short) {\it with $\bR$-polarizable graded quotients on $S$} is $(H_\bQ,W,H_\cO,F)$ consisting of locally constant sheaf $H_\bQ$ with an increasing filtration $W$ of $H_\bQ$ on $(S^{\log},\cO_S^{\log})$, locally free sheaf $H_\cO$ with a decreasing filtration $F$ of $H_\cO$ on $(S,\cO_S)$ such that $\gr_F^p$ is locally free for all $p$, 
and a specified isomorphism $\cO_S^{\log}\otimes_\bQ H_{\bQ}\simeq\cO_S^{\log}\otimes_{\cO_S}H_\cO$, whose pullbacks to each fs log point $s\in S$ satisfy the following conditions (1)--(3).

(1)  (admissibility) Let $N_1,\dots,N_n$ be a generator of the local monodromy cone $C(s):=\Hom(M_s/\cO_s^\times, \bR_{\geq 0})\subset\pi_1(s^{\log})$.
The relative monodromy weight filtrations $M(N_j, W)$ exists for all $1\le j\le n$.

(2) (Griffiths transversality) $
\nabla F^p \sub \omega_s^{1,\log}\otimes_{\cO_s}F^{p-1} \quad 
\text{for all $p\in\bZ$}.$

(3) ($\bR$-polarizablility on graded quotients) 
For each $w\in\bZ$, there is a non-degenerate $(-1)^w$-symmetric bilinear form $\langle\;,\,\rangle_w:H(\gr^W_w)_\bR\times H(\gr^W_w)_\bR\to\bR$ over $\bR$ such that the quadruple $(H(\gr^W_w), \langle\;,\,\rangle_w, H(\gr^W_w)_\cO, F(\gr^W_w))$ is an $\bR$-polarized log Hodge structure of weight $w$ on $s$.
The last part means as follows.
Let $q_1,\dots,q_r\in M_s\smallsetminus\cO_s^\times$ whose classes generate the monoid $M_s/\cO_s^\times$. 
For $t\in s^{\log}$ and $a\in\text{sp}(t)$ with $\exp(a(\log(q_j)))$ sufficiently small for all $1\le j\le r$, $(H(\gr^W_w), \langle\;,\,\rangle_w, H(\gr^W_w)_\cO, F(\gr^W_w)(a))$ is an $\bR$-polarized Hodge structure.
Here we use $H(\gr^W_w)_\bR:=\bR\otimes_\bQ H(\gr^W_w)$, $H(\gr^W_w)_\cO:=\cO_s\otimes_\bQ H(\gr^W_w)$, $F(\gr^W_w):=F(H(\gr^W_w)_{\cO})$.
\medskip

Note that, in \cite{KU09} Definition 2.4.7, \cite{KNU} Part III 1.3.2, rational polarizations 
on graded quotients were used.
But, in the present paper, we use $\bR$-polarizablility on graded quotients.
Even under this latter condition, the proof of \cite{KNU} Part III  Proposition 1.3.4 works.

\end{para}

\begin{para}\label{LMH1}

{\bf Definition.}
Given $(G,h_0)$ and $\Gamma$ as in \ref{Gamma1}.
Let $S\in\cB(\log)$.

A {\it $G$-log mixed Hodge structure with a $\Gamma$-level structure on $S$} is $(H,\mu)$ consisting of an exact $\otimes$-functor $H:\Rep(G)\to\text{LMH}(S); (V,W)\mapsto(V,W,F)$ and a global section $\mu$ of the quotient sheaf $\Gamma\backslash {\Cal I}$.

\medskip

Here $\Cal I$ is the following sheaf on $S^{\log}$. 
For an open set $U$ of $S^{\log}$, ${\Cal I}(U)$ is the set of all isomorphisms $H_\Q|_U\overset{\sim}\to \text{id}$ of $\otimes$-functors from $\Rep(G)$  to the category of local systems of $\bQ$-modules $V$ on $U$ preserving the weight filtration $W$.
\medskip

\end{para}

\begin{para}\label{Stype}
Let $(G,h_0)$ be as in \ref{Gamma1} and let $\Gamma$ and $\Sigma$ be as in \ref{sc}.

  A $G$-LMH $H$ on $S$ with a $\Gamma$-level structure $\mu$  is said to be {\it of type $(h_0,\Sigma)$} if the following (i) and (ii) are satisfied for any $s\in S$ and any $t\in s^{\log}$.
  Take a
 $\otimes$-isomorphism $\tilde \mu_t: H_{\Q,t}\overset\sim\to  \;\text{id}$ which belongs to $\mu_t$.

  (i) There is a $\sig\in \Sig$ such that the logarithm of the action of the cone $\Hom((M_S/\cO^\times_S)_s, \N)\subset \pi_1(s^{\log})$ on $H_{\Q,t}$ is contained, via $\tilde \mu_t$, in $\sigma \subset \Lie(G)_\R$.

  (ii) Let $\sig\in \Sig$ be the smallest cone satisfying (i). Let  $a: \cO_{S,t}^{\log}\to \C$ be a ring homomorphism which induces the evaluation $\cO_{S,s}\to \C$ at $s$ and consider the Hodge filtration $F$ of the functor $V\mapsto {\tilde \mu}_t a(H(V))$ in $Y$. 
Then this functor belongs to $\Dc$ and $(\sig, F)$ generates a nilpotent orbit.
\medskip
  
  If $(H, \mu)$ is of type $(h_0,\Sigma)$, we have a map $S \to \Gamma \bs D_{\Sig}$, called the {\it period map} associated to $(H, \mu)$, which sends $s\in S$ to the class of the nilpotent orbit $(\sig, Z)\in D_{\Sig}$ which is obtained in (ii).

\end{para}

\begin{para}\label{strong}
Let $(G,h_0)$ be as in \ref{Gamma1} and let $\Gamma$ and $\Sigma$ be as in \ref{sc}.

Introduce on $\Gamma\bs D_\Sigma$ the {\it strong topology}, that is the strongest topology for which the period map $S\to\Gamma\bs D_\Sigma$ is continuous for all $(S,H,\mu)$, and introduce a sheaf of holomorphic functions $\Cal O$ and a log structure $M$.

\end{para}

\begin{thm}\label{thm1}
Let $(G, h_0, \Gamma, \Sig)$ be as in \ref{strong}. Assume that $h_0$ is $\R$-polarizable (\ref{pol}).
Then

$(1)$ $\Gamma \bs D_{\Sig}$ is Hausdorff.
\smallskip

From hereafter, assume that $\Gamma$ is neat.

$(2)$ $\Gamma \bs D_{\Sig}$ is a log manifold ({\rm \cite{KNU} Part III 1.1.5}). 
In particular, $\Gamma \bs D_{\Sig}$ belongs to $\cB(\log)$.

$(3)$ $\Gamma \bs D_{\Sig}$ represents the contravariant functor from $\cB(\log)$ to {\rm (Set):}

$S\mapsto \{$isomorphism class of 
$G$-LMH over $S$ with a $\Gamma$-level structure of type $(h_0,\Sigma)$ $\}$.

$(4)$ Let $S$ be a connected, log smooth, fs log analytic space, and let $U$ be the open subspace of $S$ consisting of all points of $S$ at which the log structure of $S$ is trivial.
Assume that $S\smallsetminus U$ is a smooth divisor.

Let $(H,\mu)$ be a $G$-MHS over $U$ of  type $h_0$ (\ref{D}) endowed with a $\Gamma$-level structure (\ref{LMH1}).
Let $\varphi: U\to \G\bs D$ be the associated period map. 
Assume that $(H,\mu)$ extends to a $G$-LMH over $S$ with a $\Gamma$-level structure of type $(h_0,\Sigma)$.

Then, $\varphi$ extends to a morphism $\overline{\varphi}$ in $\cB(\log)$ in the following commutative diagram:  
$$
\CD
S@>{\overline{\varphi}}>>\G \bs D_\Sig\\
\bigcup&&\bigcup\\
U@>{\varphi}>>\G \bs D.
\endCD
$$
\end{thm}

\section{Log higher Albanese manifolds}\label{s:lha}
\medskip

We review here formulations and results of higher Albanese manifolds in \cite{HZ87} and of log higher Albanese manifolds in \cite{KNU17}.

\begin{para}\label{Lie}

Let $X$ be a connected smooth algebraic variety over $\bC$.
Fix $b\in X$. 
Let $\Gamma$ be a quotient group of $\pi_1(X, b)$ which is  torsion-free and nilpotent. 
\smallskip

Let $\cG=\cG_\Gamma$ be the unipotent algebraic group over $\Q$ whose Lie algebra is defined as follows:
Let $I$ be the augmentation ideal $\text{Ker}(\bQ[\Gamma]\to \bQ)$ of $\Q[\Gamma]$. 
Then $\Lie(\cG)$ is the $\Q$-subspace of $\Q[\Gamma]^{\wedge}:=\varprojlim_n \Q[\Gamma]/I^n$ generated by all 
$\log(\gamma)$ ($\gamma \in \Gamma$). The Lie product of $\Lie(\cG)$ 
is defined by $[x,y]= xy-yx$. 
We have $\Gamma \subset \cG(\Q)$. 

\end{para}

\begin{para}\label{hAlb}

Let $\pi_1=\pi_1(X, b)$.
Let $J$ be the augmentation ideal $\text{Ker}(\bQ[\pi_1]\to \bQ)$ of $\Q[\pi_1]$.
For a positive integer $n$, let $\Gamma_n$ be the image of $\pi_1\to\bQ[\pi_1]/J^n$.
Then $\Lie(\cG_{\Gamma_n})$ has a mixed Hodge structure induced from de Rham theory on the path spaces over $X$ by Chen's iterated integral.

For a given $\Gamma$ as in \ref{Lie}, there exists $n\ge1$ such that $\Gamma$ is a quotient of 
$\Gamma_n$.
Hereafter we assume that $\Lie(\cG_{\Gamma})$ has a quotient mixed Hodge structure 
of the one on $\Lie(\cG_{\Gamma_n})$.
Note that this mixed Hodge structure on $\Lie(\cG_{\Gamma})$ is independent of the choice of $n$. 

We note that there is an insufficient statement on mixed Hodge structre on $\Lie(\cG_{\Gamma})$ in \cite{KNU17} 6.1.2.
The authors of \cite{KNU17} agreed to correct this part, so as to assume the existence of this mixed Hodge structure on $\Lie(\cG_{\Gamma})$ as above in the present paper.

Let $\cG=\cG_\Gamma$.
Let $F^0\Lie(\cG)_\bC$ be the $0$-th Hodge filter on $\Lie(\cG)_\bC$ and let $F^0\cG(\bC)$ be the corresponding subgroup of $\cG(\C)$.
The {\it higher Albanese manifold} associated to $(X,\Gamma)$ is defined in \cite{HZ87} as 
$$
A_{X,\Gamma}:= \Gamma\bs \cG(\C)/F^0\cG(\C).
$$

\end{para}

\begin{para}\label{G}

Take a $\Q$-MHS $V_0$ with polarizable $\gr^{W}$ having the $\Q$-MHS on $\Lie(\cG)_\bQ$ 
with $\cG=\cG_\Gamma$ in \ref{hAlb} as a direct summand.

Let $Q\subset \Aut(V_{0,\Q})$ be the {\it Mumford-Tate group} associated to $V_0$, i.e., the Tannaka group of the Tannaka category $\langle V_0 \rangle$ generated by $V_0$: $\langle V_0 \rangle \overset\sim\to \Rep(Q)$.
Explicitly, it is the smallest $\Q$-subgroup $Q$ of $\Aut(V_{0,\Q})$ such that $Q_\R$ contains the image of the homomorphism $h:S_{\C/\R} \to \Aut(V_{0,\R})$ and such that $\Lie(Q)_\R$ contains $\delta$. 
Here $h$ and $\delta$ are determined by the canonical splitting of the $\bQ$-MHS $V_0$ (\cite{CKS87}, \cite{KNU} Part II 1.2).
\bigskip  

The action of $Q$ on $\Lie(\cG)$ induces an action of $Q$ on $\cG$. 
By this, define a semi-direct product $G$ of $Q$ and $\cG$:
$$
1\to \cG\to G\to Q\to 1.
$$
We have $\cG\subset G_u$. 
We have $h_0: S_{\C/\R}\to (Q_{\red})_\R= (G_{\red})_\R$ by using the Hodge decomposition of $\gr^WV_0$.

\end{para}

\begin{para}\label{moduli}

Let $D_G$ (resp.\ $D_Q$) be the period domain $D$ for $G$ (resp.\ $Q$) and $h_0$ in \ref{G}. 
We have a canonical map $\Gamma \bs D_G\to D_Q$ induced by $G\to Q$. 

Let $b_Q\in D_Q$ be  the isomorphism class of the evident functor $\Rep(Q)\to\text{$\Q$-MHS}$ by  definition in \ref{G}, and let $b_G\in D_G$ be the isomorphism class of $\Rep(G)\to\Rep(Q)\overset{b_Q}\to\text{$\Q$-MHS}$ via the section $Q\hookrightarrow G$.

Let $\cD$ be the fiber of the map $D_{G} \to D_{Q}$ over $b_Q$.

\end{para}

\begin{thm}\label{t:quot}
The map $G_u(\C) \to D_G\; ;\; g \mapsto g\cdot b_G$ 
induces an isomorphism 
$A_{X,\Gamma}=\Gamma\bs \cG(\C)/F^0\cG(\C)\simeq\Gamma\bs \cD$ 
of analytic manifolds.
\end{thm}

\begin{para}\label{HZ}

Let $\cC_{X,\Gamma}$ be the category of variations of $\Q$-MHS $\cH$ on $X$ satisfying the following three conditions:

(1) For any $w\in\bZ$, $\gr^W_w\cH$ is a constant polarizable Hodge structure.

(2) $\cH$ is {\it good at infinity} in the sense of \cite{HZ87} (1.5), i.e., there exists a smooth compactification $\overline{X}$ of $X$ with normal crossing boundary divisor $\overline{X}\smallsetminus X$ such that
the Hodge filtration bundles extend to sub-bundles of the canonical extension of $\cO$-module of 
$\cH$ which induce the corresponding thing for each $\gr^W_w\cH$, and that, for the nilpotent logarithm $N_j$ of a local monodromy transformation about a component of $\overline{X}\smallsetminus X$, the relative monodromy weight filtration $M(N_j,W)$ exists.

(3) The monodromy action of $\pi_1(X,b)$ factors through $\Gamma$.

\end{para}

Hain and Zucker showed 

\begin{thm}\label{t:HZ}{\rm (\cite{HZ87} (1.6) Theorem).}
The category $\cC_{X,\Gamma}$ is equivalent to the category of $\bQ$-MHS $V$ with polarizable $\gr^WV$ endowed with an action of $\Lie(\cG)$ such that $\Lie(\cG)\otimes V\to V$ is a homomorphism of MHS.
\end{thm}

\begin{para}\label{FG}

Define a contravariant functor $\cF_{\Gamma}:\cB(\log) \to\text{Sets}$ as follows: 
For $S\in\cB(\log)$, $\cF_{\Gamma}(S)$ is the set of isomorphism classes of pairs $(H,\mu)$ of an exact $\otimes$-functor $H:\cC_{X,\Gamma} \to \text{MHS}(S)$ and a $\Gamma$-level structure $\mu$ satisfying the following condition (i).
Here a {\it $\Gamma$-level structure} means a global section 
of the sheaf $\Gamma \bs \Cal I$, where $\Cal I$ is the sheaf of functorial 
$\otimes$-isomorphisms $H(\cH)_{\Q} \overset{\sim}\to \cH(b)_{\Q}$ of $\Q$-local systems preserving weight filtrations.

 (i) For any $\Q$-MHS $h$, we have a functorial $\otimes$-isomorphism $H(h_X) \cong h_S$ such that the induced isomorphism of local systems $H(h_X)_\Q\cong h_\Q=h_X(b)_\Q$ belongs to $\mu$. Here $h_X$ (resp.\  $h_S$) denotes the constant variation (resp.\ family) of $\Q$-MHS over $X$ (resp.\ $S$) associated to $h$.

\begin{thm}\label{t:FG}
Let the notation be as in \ref{FG}.
The functor $\cF_{\Gamma}$ is represented by $A_{X,\Gamma}\simeq\Gamma\bs\cD$. 
\end{thm}

This follows from Theorem \ref{t:quot} and Theorem \ref{t:HZ}.

Let $\varphi:X\to A_{X,\Gamma}$ be the higher Albanese map.

\end{para}

\begin{para}\label{FGS}

Let $\Sig$ be a weak fan in $\Lie(G)$ such that $\sig\subset \Lie(\cG)_\R$ for any $\sig\in \Sig$. 
Assume that $\Sigma$ and $\Gamma$ are strongly compatible.
Let $\Gamma \bs D_{G,\Sig}\to D_Q$ be a canonical morphism induced by $G\to Q$.
Define 
$$
\text{$A_{X, \Gamma,\Sig}$:= (the fiber of $\Gamma\bs D_{G,\Sigma}\to D_Q$ over $b_Q$) $\in\cB(\log)$}
$$

Define a contravariant functor $\cF_{\Gamma,\Sig}:\cB(\log) \to\text{Sets}$ as follows:
For $S\in\cB(\log)$, ${\cF}_{\Gamma,\Sig}(S)$ is the set of isomorphism classes of pairs $(H,\mu)$ 
consisting of an exact $\otimes$-functor $H:\cC_{X,\Gamma}\to\LMH(S)$ and  a $\Gamma$-level structure $\mu$ satisfying the condition (i) in \ref{FG} and also the following condition (ii).
 
(ii) The following (ii-1) and (ii-2) are satisfied for any $s\in S$ and any $t\in s^{\log}$. 
Let $\tilde \mu_t: H(\cH)_{\Q,t}\cong \cH(b)_\Q$ be a  functorial $\otimes$-isomorphism which belongs to $\mu_t$.

  (ii-1) There is a $\sig\in \Sig$ such that the logarithm of the action of the local monodromy cone 
$\Hom((M_S/\cO^\times_S)_s, \N)\subset \pi_1(s^{\log})$ on $H_{\Q,t}$ is contained, via $\tilde \mu_t$, in $\sigma \subset \Lie(\cG)_\R$.

  (ii-2) Let $\sig\in \Sig$ be the smallest cone which satisfies (ii-1) and let $a: \cO_{S,t}^{\log}\to \C$ be a ring homomorphism which induces the evaluation $\cO_{S,s}\to \C$ at $s$. Then, for each $\cH\in \cC_{X,\Gamma}$, $(\sig, \tilde \mu_t(a(H(\cH))))$ generates a nilpotent orbit in the sense of \cite{KNU} Part III, 2.2.2.

\begin{thm}\label{t:FGS}
Let the notation be as in \ref{t:FG} and \ref{FGS}.

$(1)$ The functor $\cF_{\Gamma,\Sig}$ is represented by $A_{X, \Gamma,\Sig}$.

$(2)$ Let $\overline{X}$ be a smooth algebraic variety over $\C$ which contains $X$ 
as a dense open subset such that the complement 
$\overline{X}\smallsetminus X$ is a smooth divisor. 
Endow $\overline{X}$ with the log structure associated to this divisor. 
Assume that $\Sigma$ is the fan consisting of all rational nilpotent cones in $\Lie(\cG)_\bR$ of rank $\le1$ (denoted by $\Xi$ in \cite{KNU17} 6.2.5).
Then, the {\rm higher Albanese map} $\varphi: X \to A_{X,\Gamma}$ extends 
uniquely to a morphism $\overline{\varphi}: \overline{X} \to 
A_{X,\Gamma,\Sigma}$ of log manifolds.

\end{thm}

Since an object of $\cC_{X, \Gamma}$ is good at infinity (\ref{HZ}), it extends to an LMH over $\overline{X}$.
Hence (2) follows from (1)  and the general theorem \ref{thm1} (4).

\end{para}

\section{Description of a result of Deligne by log higher Albanese map}\label{s:example}

For a group $\Gamma^{(n)}$ in \ref{Gamma^n} below, Deligne \cite{D89} showed that polylogarithms appear in the higher Albanese map $X\to A_{X,\Gamma^{(n)}}$ (cf.\ Section A below).
Here we describe them in our framework in \cite{KNU17} (Section 2 in the present paper).

\begin{para}\label{bsetting}
Let $X:=\bP^1(\bC)\smallsetminus\{0,1,\infty\}\subset\overline{X}:=\bP^1(\bC)$ with affine coordinate $x$.
Let $b:=(0,1)$ the \lq\lq tangential base point" over $0\in\overline{X}$ with tangent $v_0\in T_0(\overline{X})=\Hom_\bC(m_0/m_0^2,\bC)$ defined by $v_0(x)=1$ in \cite{D89} Section 15.
This is understood in log geometry in the following way.
Let $y=(0,h)\in\overline{X}^{\log}$ be the point lying over $0\in X$, where $h:M^{\gp}_{\overline{X},0}=\cO_{\overline{X},0}^{\times}\times x^\bZ \to {\bf S}^1$ is the argument function which is a group homomorphism sending $f\in\cO_{\overline{X},0}^{\times}$ to $f(0)/|f(0)|$ and $x$ to $v_0(x)/|v_0(x)|=1$ (\cite{KNU17} 6.3.7).
Let $u_0\in\cO_{\overline{X},y}^{\log}$ be the branch of $\log(x)$ having real value on $\bR_{>0}$.
(This $u_0$ is the branch denoted by $f\in\cO_{\overline{X},y}^{\log}$ in \cite{KNU17} 6.3.7 (ii), and $u_0$ can be also regarded as the function $2\pi iz$ on $\tilde S^{\log}$ in 1.1.1.)
Then the corresponding base point in the boundary in our sense is  $b=(y,a)$, where $a:\cO^{\log}_{\overline{X},y}=\bC\{x\}[u_0]\to \C$ is  the specialization which is a ring homomorphism sending $x$ to $0$ and $u_0$ to $a(u_0)=\log(v_0(x))=\log(1)=0$ (\cite{KNU17} 6.3.7 (ii)).
\medskip

See \cite{KNU17} 6.3.6, 6.3.7 for more general description of the above  correspondence of boundary points.

\end{para}

\begin{para}\label{Gamma}

The inclusion $X\subset\bG_m(\bC)=\bC^{\times}$ induces $\pi_1(X,b)\to\pi_1(\bG_m(\bC),b)=\bZ(1)$.
Let $K$ be its kernel, and let $\Gamma:=\pi_1(X,b)/[K,K]$ and  
$\Gamma_1:=K/[K,K]$.
Then
$$1\to\Gamma_1\to\Gamma\to\bZ(1)\to1.$$

\end{para}

\begin{para}\label{Gamma^n}

Let $Z^n\Gamma$ be the descending central series of $\Gamma$ defined by  $Z^{n+1}\Gamma:=[Z^n\Gamma,\Gamma]$ starting with $Z^1\Gamma=\Gamma$.

Let $\Gamma^{(n)}:=\Gamma/Z^{n+1}(\Gamma)$ and
$\Gamma_1^{(n)}:=\Image\,(\Gamma_1\to\Gamma^{(n)})$.
Let $\gamma_0,\gamma_1\in\Gamma^{(n)}$ be the classes of small loops anticlockwise around $0$ and clockwise around $1$, respectively.
Then, we have 
\vskip5pt

$\Gamma^{(n)}=\langle \gamma_{0},\gamma_{1}\rangle$,
$({\rm ad}\,\gamma_{0})^{k-1}\gamma_{1}$  $(1\le k\le n)$ are commutative, 
$\Gamma_{1}^{(n)}=\sum_{k=1}^n\bZ({\rm ad}\,\gamma_{0})^{k-1}\gamma_{1}$.

\end{para}

\begin{para}\label{DLam}

Let $\Lambda=(V,W,(\lan\;,\ran_w)_{w\in\bZ},(h^{p,q})_{p,q\in\bZ})$ be as follows.
$V$ is a free $\bZ$-module with basis $e_1,e_2,e_3,\dots,e_{n+1}$.
$W$ is a weight filtration on $V_\bQ$ defined by
$$
W_{-2n-1}=0\subset W_{-2n}=W_{-2n+1}=\bQ e_1
\subset W_{-2n+2}=W_{-2n+3}=W_{-2n+1}+\bQ e_2
$$
$$
\subset\cdots\subset W_0=W_{-1}+\bQ e_{n+1}=V_{\bQ}.
$$
$\langle\;,\,\rangle_w : \gr^W_w(V_\bQ) \times \gr^W_w(V_\bQ)\to \Q$  ($w\in\bZ$) are the $\Q$-bilinear forms characterized by $\langle e_{n+1+k},e_{n+1+k}\rangle_{2k}=1$ for $k=0,-1,\dots,-n$.
$h^{k,k}=1$ for $k=0,-1,\dots,-n$, and $h^{p,q}=0$ for the other $(p,q)$. 

Let $D(\Lambda)$ be the period domain in \cite{KNU} Part III with universal Hodge filtration $F$:
$$
F^{1}=0\subset F^0=\bC(e_{n+1}+\sum_{n\ge j\ge1}a_{j,n+1}e_j)
\subset
F^{-1}=F^0+\bC(e_n+\sum_{n-1\ge j\ge1}a_{j,n}e_j) 
$$
$$
\subset\cdots\subset F^{-n}=F^{-n+1}+\bC e_1=V_{\bC}.
$$

\end{para}

\begin{para}\label{action}

Let $\cG$ be the unipotent group $\cG$ in \ref{Lie} for $\Gamma^{(n)}$.
Define an action of $\Lie(\cG)$ on $V_\bQ$ by $N_0=\log(\gamma_0)$, $N_1=\log(\gamma_1)$:
$$
N_0e_j=e_{j-1}\; (j=2,\dots,n),\;\;N_0e_j=0\;(j=1,n+1),
$$
$$
N_1e_{n+1}=-e_n,\;\; N_1e_j=0\;(j=1,2,\dots,n).
$$

Then
$$
(-N_0+N_1)^j=(-N_0)^j+(-\Ad N_0)^{j-1}N_1\quad
(1\le j\le n+1).
$$

\end{para}

\begin{para}\label{hAlb2}

Let $X$ be as in \ref{bsetting} and $\Gamma^{(n)}$ be as in \ref{Gamma^n}.
We consider the higher Albanese manifold $A_{X,\Gamma^{(n)}}$ of $X$ by using the base point $b$ in \ref{bsetting}. 

The $\bQ$-MHS on $\Lie(\cG)$ is as follows: 
$N_0$ and $N_1$ are of Hodge type $(-1,-1)$ and compatible with bracket and hence $F^0\cG(\bC)=\{1\}$.
Thus the higher Albanese manifold is
$$
A_{X,\Gamma^{(n)}}=\Gamma^{(n)}\bs\cG(\bC).
$$

\end{para}

\begin{lem}\label{l:MHSact}
Let $F$ and $N_j$ $(j=0,1)$ be as in \ref{DLam} and in \ref{action}.

\noindent
{\rm(i)} We have the following.

$(1)$ $(N_0,F)$ satisfies the Griffiths transversality if and only if 
$$
a_{k,n+1}=0 \quad (2\le k\le n);\quad
a_{1,k}=a_{l-k+1,l}\quad(2\le k< l\le n). 
$$

$(2)$ $(N_1,F)$ satisfies the Griffiths transversality if and only if 
$$
a_{k,n}=0 \quad (1\le k\le n-1).
$$

$(3)$ $(-N_0+N_1,F)$ satisfies the Griffiths transversality if and only if 
$$
a_{1,k}=a_{l-k+1,l}\quad(2\le k<l\le n+1). 
$$

\noindent
{\rm(ii)} The following three conditions are equivalent.

$(1)$ The Lie action $\Lie(\cG)\otimes V\to V$ in \ref{action} is a homomorphism of MHS with respect to the MHS on $\Lie(\cG)$ in \ref{hAlb2} and the MHS $(V,W,F)$ in \ref{DLam}.

$(2)$ For $j=0$ and $1$, $(N_j,F)$ satisfies the Griffiths transversality.

$(3)$ $a_{j,k}=0$ unless $(j,k)=(1,n+1)$.

\end{lem}

The assertions are easily verified by direct computation.

\begin{para}\label{hAlb3}

For any fixed $a\in\bC$, denote by $F(a)$ the Hodge filtration in \ref{l:MHSact} (ii) (3) with $a_{1,n+1}=a$.
By the action in \ref{action}, we define 
$$
\cD:=\exp\left(\bC N_0+\sum_{k=1}^n\bC(\Ad N_0)^{k-1}N_1\right)F(a)\subset D(\Lambda).
$$
Then, this $\cD$ coincides with $\cD$ in \ref{moduli}.
Hence $\cG(\bC)\simeq\cD$ and $A_{X,\Gamma^{(n)}}\simeq\Gamma^{(n)}\bs\cD$ as complex analytic manifolds.

\end{para}

\begin{para}\label{eq}

Let  $\varphi:X\to A_{X,\Gamma^{(n)}}\simeq\Gamma^{(n)}\bs\cD$ be the composite of higher Albanese map and the isomorphism in \ref{hAlb3}.
Let $F(x)$ be the pullback by $\varphi$ of the universal Hodge filtration on $\Gamma^{(n)}\bs\cD$.

Since $F(x)$ is rigid by Theorem \ref{t:HZ}, we consider a connection equation:
$$
dF(x)=\omega F(x),\quad
\omega:=(2\pi i)^{-1}\frac{dx}{x}N_0+(2\pi i)^{-1}\frac{dx}{1-x}N_1.
$$
That is,
$$
da_{k-1,k}(x)=(2\pi i)^{-1}\frac{dx}{x}\quad
(2\le k\le n),
$$
$$
da_{n,n+1}(x)=-(2\pi i)^{-1}\frac{dx}{1-x},
$$
$$
da_{j,k}(x)=(2\pi i)^{-1}a_{j+1,k}(x)\frac{dx}{x}\quad
(3\le k\le n+1,\;1\le j\le k-2).
$$

\end{para}

\begin{para}\label{sol}

This system is solved by iterated integrals.
The solutions are
$$
a_{j,k}(x)=\frac{1}{(k-j)!}((2\pi i)^{-1}\log (x))^{k-j}\quad
(2\le k\le n,\;1\le j\le k-1),
$$
$$
a_{j,n+1}(x)=-(2\pi i)^{-n-1+j}l_{n+1-j}(x)\quad(1\le j\le n).
$$
Here the $l_j(x)$ are polylogarithms, in particular $l_1(x)=-\log(1-x)$.

\end{para}

\begin{para}\label{table}

Table of solutions:
$$
\left(\CD
1\;\;&a_{1,2}\;\;&\ldots&a_{1,n}&a_{1,n+1}\\
0&1&\ddots&\vdots&\vdots\\
\vdots&0&\ddots\;\;&\;a_{n-1,n}\;\;&a_{n-1,n+1}\\
\vdots&\vdots&\ddots&1&a_{n,n+1}\\
0&0&\ldots&0&1
\endCD\right)=
\left(\CD
1\;\;&(2\pi i)^{-1}\log(x)\;&\ldots\;\;&\frac{((2\pi i)^{-1}\log(x))^{n-1}}{(n-1)!}\;\;&-(2\pi i)^{-n}l_n(x)\\
0&1&\ddots&\vdots&\vdots\\
\vdots&0&\ddots&(2\pi i)^{-1}\log(x)&-(2\pi i)^{-2}l_2(x)\\
\vdots&\vdots&\ddots&1&-(2\pi i)^{-1}l_{1}(x)\\
0&0&\ldots&0&1
\endCD\right).
$$

Note that, for $1\le j\le n$, 
$$
\exp((2\pi i)^{-1}\log(x)N_0)e_j
=e_j+(2\pi i)^{-1}\log(x)e_{j-1}+\cdots
+\frac{1}{(j-1)!}((2\pi i)^{-1}\log (x))^{j-1}e_1,
$$
for $j=n+1$,
$$
\exp\left(-\sum_{n\ge k\ge1}(2\pi i)^{-k}l_{k}(x)(\Ad N_0)^{k-1}N_1\right)
e_{n+1}
=e_{n+1}-\left(\sum_{n\ge k\ge1}(2\pi i)^{-k}l_{k}(x)e_{n+1-k}\right).
$$

\end{para}

\begin{para}\label{Hfilt}
For $\alpha,\beta,\lambda_{2},\dots,\lambda_{n}\in\bC$, let 
$F=F(\alpha,\beta,\lambda_{2},\dots,\lambda_{n})$ 
be a Hodge filtration:
$$
F^{1}=0\subset F^0=\bC(e_{n+1}+\beta e_{n}+
\lambda_{2}e_{n-1}+\dots+\lambda_{n}e_{1})
$$
$$
\subset
F^{-1}=F^0+\bC\left(e_n+\alpha e_{n-1}+\frac{\alpha^{2}}{2!}e_{n-2}
+\dots+\frac{\alpha^{n-1}}{(n-1)!}e_1\right)\subset\cdots
$$
$$
\subset F^{-n+1}=F^{-n+2}+\bC(e_2+ \alpha e_1)
\subset F^{-n}=F^{-n+1}+\bC e_1=V_{\bC}.
$$

\end{para}

\begin{para}\label{hAmap}

Let $\varphi:X\to A_{X,\Gamma^{(n)}}\simeq \Gamma^{(n)}\bs\cD$ be the higher Albanese map in \ref{eq}.
We have a commutative diagram
$$
\CD
&&\tilde\varphi(X)&\;&\subset& \cD\\
&\tilde\varphi\nearrow& @VVV&@VVV\\
X@>\sim>\varphi>\varphi(X)&\;&\subset \;&\Gamma^{(n)}\bs\cD
\endCD
$$
where $\tilde\varphi:X\to\cD$ is a multi-valued map corresponding to the Hodge filtration
$$
x\mapsto F((2\pi i)^{-1}\log (x),-(2\pi i)^{-1}l_1(x),\dots, -(2\pi i)^{-n}l_{n}(x))
$$
in the nottion in \ref{Hfilt}.
$\tilde\varphi(X)\to X$ and $\tilde\varphi(X)\to\varphi(X)$ are $\Gamma^{(n)}$-torsors and $\varphi:X\overset\sim\to\varphi(X)$ is an isomorphism.

\end{para}

\begin{para}\label{loghA}

Let $\Sigma$ be the set of all cones of the form $\R_{\geq 0} N$ with $N\in\Lie(\cG)$. 
We consider the extended period domain $D(\Lam)_{\Sigma}$ in \cite{KNU} Part III. 
This is only a set.
By using the strong topology (\cite{KU09} Section 3.1), 
the quotient $\Gamma^{(n)}\bs D(\Lambda)_{\Sigma}$ has a structure of a log manifold. 
Define $\Gamma^{(n)}\bs\cD_\Sigma$ to be the closure of $\Gamma^{(n)}\bs\cD$ in $\Gamma^{(n)}\bs D(\Lambda)_\Sigma$.
This inherits a structure of log manifold.
We have $A_{X,\Gamma^{(n)},\Sigma}\simeq\Gamma^{(n)}\bs\cD_\Sigma$ in the category $\cB(\log)$.

Let $N\in\Lie(\cG)$ and $\sigma:=\R_{\geq 0} N$.
Let $\Gamma_{\sigma}$ be the group generated by the monoid $\Gamma^{(n)}\cap\exp(\sigma)$.
If we use as $\Sigma$ the fan consisting of the cone $\sigma$ and $0$, also denoted by $\sigma$ by abuse of notation, we have $A_{X,\Gamma_{\sigma},\sigma}\simeq\Gamma_{\sigma}\bs\cD_\sigma$ in the category $\cB(\log)$.

\end{para}

\begin{para}\label{local0}
Let $N_0$ be as in \ref{action} and set $\sig_0=\R_{\geq 0}N_0$. 
Let $F=F(\alpha,\beta,\lambda_{2},\dots,\lambda_{n})$ be as in \ref{Hfilt}.
By Lemma \ref{l:MHSact} (i) (1), $(N_0,F)$ satisfies the Griffiths transversality if and only if $\beta=\lambda_{2}=\dots=\lambda_{n-1}=0$.
If this is the case, $(N_0,F)$ generates a $\sigma_0$-nilpotent orbit, since admissibility and $\bR$-polarizability on $\gr^W$ trivially hold.
We describe the local structure of $\Gamma_{\sigma_0}\bs\cD_{\sigma_0}$ near the image $p_0$ of this nilpotent orbit.

Let $Y:=\{(q, \beta, \lambda_2,\dots,\lambda_n)\in \C^{n+1}\;|\; \beta=\lambda_{2}=\dots=\lambda_{n-1}=0 \;\text{if}\; q=0\}$ 
be the log manifold with the strong topology, with the structure sheaf of rings which is the inverse image of the sheaf of holomorphic functions on $\C^{n+1}$, and with the log structure generated by $q$. 
Then there is an open neighborhood $U$ of $(0,0,\dots,0,\lambda_n)$ in $\C^{n+1}$ and an open immersion 
$$
Y\cap U \hookrightarrow \Gamma_{\sigma_0}\bs\cD_{\sigma_0}
$$ 
of log manifolds which sends $(q, \beta, \lambda_2,\dots,\lambda_n)\in Y\cap U$ with $q\neq 0$ to the class of $F(\alpha,\beta,\lambda_2,\dots,\lambda_n)$,  
where $\alpha\in \C$ is such that $q=e^{2\pi i\alpha}$, and which sends $(0,0,\dots,0,\lambda_n)$ to $p_0$. 

\end{para}

\begin{para}\label{naive0}

Near $x=0$, a nilpotent orbit in naive sense is
$$
(1)\hskip40pt
\exp((2\pi i)^{-1}\log(x)N_0)F(0,0,\dots,0,\lambda_n^0)
=F((2\pi i)^{-1}\log(x),0,\dots,0,\lambda_n^0),\hskip40pt
$$
where $\lambda_n^0=-(2\pi i)^{-n}l_n(0)$.
The corresponding \lq\lq higher Albanese map" (i.e., local version about $0$ of $\tilde\varphi$ in \ref{hAmap}) is
$$
(2)\hskip90pt
F((2\pi i)^{-1}\log(x),-(2\pi i)^{-1}l_1(x),\dots,-(2\pi i)^{-n}l_n(x))
\hskip90pt
$$
under the condition $l_j(0)=0$ ($1\le j\le n-1$).
These two are asymptotic when $x$ goes to the boundary point $b=(y,a)$ with $y=(0,h)\in\overline{X}^{\log}$ and $a$ being the specialization at $y$ in \ref{bsetting}.

\end{para}

\begin{para}\label{over0}

As above, let $u_0$ be the branch of $\log (x)$ in \ref{bsetting} and $T$ an indeterminate over $\cO_{\overline{X},0}$.
Then, by 1.1.1, we have an isomorphism $\cO_{\overline{X},y}^{\log}=\cO_{\overline{X},0}[u_0]\simeq \cO_{\overline{X},0}[T]$ of $\cO_{\overline{X},0}$-algebras under $(2\pi i)^{-1}u_0\leftrightarrow T$.
Consider an $\cO_{\overline{X},0}$-algebra homomorphism $\cO_{\overline{X},0}[T]\to\cO_{\overline{X},0}$, $T\mapsto x$.

Under the initial condition in \ref{naive0} given by the base point $b$ in \ref{bsetting}, we have
$$
l_j(x)=\sum_{k=1}^\infty \frac{x^k}{k^j}\quad(1\le j\le n-1),\quad
l_n(x)=c+\sum_{k=1}^\infty \frac{x^k}{k^n}
$$ 
on a simply connedcted neighborhood $\overline{X}_0$ of $0\in\overline{X}$, where $c:=-(2\pi i)^n\lambda_n^0$.

Let $\alpha=(2\pi i)^{-1}\log(x)$.
Then, as $x\to0$, $\exp(-\alpha N_0)(\text{$F$ in \ref{naive0} (2)})$ converges to $F(0,0,\dots,0,\lambda_n^0)$ in $\cD$ (\ref{hAlb3}), and hence the class of ($F$ in \ref{naive0} (2)) converges to the class $p_0$ (\ref{local0}) of the nilpotent orbit $(\sigma_0,\exp(\sigma_{0,\bC})F(0,0,\dots,0,\lambda_n^0))$ in 
$\Gamma_{\sigma_0}\bs\cD_{\sigma_0}$.
We thus have an extension of the higher Albanese map over $\overline{X}_0$ (Theorem \ref{t:FGS} (2)):
$$
\overline\varphi_0:\overline{X}_0\to\Gamma_{\sigma_0}\bs\cD_{\sigma_0}.
$$

This is a morphism in the category $\cB(\log)$.
The log structure on the source (resp.\ the target) is given by $x$ (resp.\ $q$).
The pullback of the universal log mixed Hodge structure on the target coincides with the log mixed Hodge structure on the source.
\end{para}

\begin{para}\label{log0}

By using log mixed Hodge theory, \ref{naive0} is described as follows.

Taking  the images of the nilpotent orbit in naive sense \ref{naive0} (1) and the \lq\lq higher Albanese map" \ref{naive0} (2), we have their real analytic extensions with boundary
$$
\overline\nu_0^{\log}, \;
\overline\varphi_0^{\log}:\overline{X}_0^{\log}\to(\Gamma_{\sigma_0}\bs\cD_{\sigma_0})^{\log}.
$$
Here, $\overline{X}_0^{\log}$ is like Example 1.1.1, and
$(\Gamma_{\sigma_0}\bs\cD_{\sigma_0})^{\log}$ coincides with the moduli of nilpotent 
$i$-orbits $\Gamma_{\sigma_0}\bs\cD_{\sigma_0}^{\sharp}$ in the present situation (\cite{KNU} III Theorem 2.5.6).

Let $\tilde{\overline{X}}_0^{\log}$ be the universal covering of $\overline{X}_0^{\log}$.
The above maps are still lifted to
$$
\tilde{\overline\nu}_0^{\log}, \;
\tilde{\overline\varphi}_0^{\log}:\tilde{\overline{X}}_0^{\log}\to\cD_{\sigma_0}^{\sharp}.
$$
The boundary point $b$ in \ref{naive0} can be understood as the point $b=(z=0+i\infty)=(u_0=-\infty+i0)\in\tilde{\overline{X}}_0^{\log}$.
We have 
$(\exp(-(2\pi i)^{-1}\log(x)N_0)(\ref{naive0} \;(2)))(b)=F(0,0,\dots,0,\lambda_n^0)$, and 
$$
\tilde{\overline\nu}_0^{\log}(b)=\tilde{\overline\varphi}_0^{\log}(b)=
(\text{nilpotent $i$-orbit generated by $(N_0, F(0,0,\dots,0,\lambda_n^0))$})\in\cD_{\sigma_0}^{\sharp}.$$

\end{para}

\begin{para}\label{local1}
Let now $\sigma_1=\bR_{\ge0}N_1$ for $N_1$ in \ref{action}.
Let $F=F(\alpha,\beta,\lambda_{2},\dots,\lambda_{n})$ be as in \ref{Hfilt}.
By Lemma \ref{l:MHSact} (i) (2), $(N_1,F)$ satisfies the Griffiths transversality if and only if $\alpha=0$.
If this is the case, $(N_1,F)$ generates a $\sigma_1$-nilpotent orbit, since admissibility and $\bR$-polarizability on $\gr^W$ trivially hold.
We have a similar description of the local structure of $\Gamma_{\sigma_1}\bs \cD_{\sigma_1}$ near the image $p_1$ of this nilpotent orbit.

Let $Y$ be the log manifold $\{(\alpha,q,\lambda_2,\dots,\lambda_n)\in \C^{n+1}\;|\; \alpha=0 \;\text{if}\; q=0\}$ with the strong topology, the structure sheaf and the log structure defined by $q$. 
Then there is an open neighborhood $U$ of $(0,0,\lambda_2,\dots,\lambda_n)$ in $\C^{n+1}$ and an open immersion 
$$
Y\cap U \hookrightarrow \Gamma_{\sigma_1}\bs \cD_{\sigma_1}
$$ 
of log manifolds which sends $(\alpha,q, \lambda_2,\dots,\lambda_n)\in Y\cap U$ with $q\neq 0$ to the class 
of $F(\alpha, \beta, \lambda_2,\dots,\lambda_n)$, where $\beta\in \C$ is such that $q=e^{2\pi i\beta}$, and which sends $(0,0,\lambda_2,\dots,\lambda_n)$ to $p_1$.

\end{para}

\begin{para}\label{naive1}
We assume the initial condition in \ref{naive0}.
Near $x=1$, a nilpotent orbit in naive sense is
$$
(1)\hskip30pt\exp((2\pi i)^{-1}\log(1-x)N_{1})\cdot
F(0,0,-(2\pi i)^{-2}\zeta(2),\dots,-(2\pi i)^{-n}(c+\zeta(n)))\hskip30pt
$$
$$
=F(0,-(2\pi i)^{-1}l_1(x),-(2\pi i)^{-2}\zeta(2),\dots,-(2\pi i)^{-n}(c+\zeta(n))).
$$
The corresponding \lq\lq higher Albanese map" (i.e., local version about $1$ of $\tilde\varphi$ in \ref{hAmap}) is
$$
(2)\hskip84ptF((2\pi i)^{-1}\log(x),-(2\pi i)^{-1}l_1(x),\dots,-(2\pi i)^{-n}l_n(x)).\hskip84pt
$$
These two are asymptotic when $x$ goes to the tangential boundary point 
$\tilde p_1:=(1,-1)$ with tangent $v_1\in T_1(\overline{X})=\Hom_\bC(m_1/m_1^2,\bC)$ defined by $v_1(1-x)=-1$.
This is the boundary point in our sense described as follows.
Let $u_1$ be the branch of $\log(1-x)$ having real value on $\bR_{<1}$.
Then the corresponding point in the boundary in our sense is  
$\tilde p_1=(y,a)$ with $y=(1,h)\in\overline{X}^{\log}$ such that 
the argument function 
$h:M^{\gp}_{\overline{X},1}=\cO_{\overline{X},1}^{\times}\times (1-x)^\bZ \to {\bf S}^1$ is a group homomorphism sending $f\in\cO_{\overline{X},1}^{\times}$ to $f(1)/|f(1)|$ and $1-x$ to $v_1(1-x)/|v_1(1-x)|=-1$, and the specialization $a:\cO^{\log}_{\overline{X},y}=\bC\{1-x\}[u_1]\to \C$ is a ring homomorphism sending $1-x$ to $0$ and $u_1$ to 
$a(u_1)=-a(-u_1)=\log(v_1(-(1-x)))=\log(1)=0$ 
(cf.\ \cite{KNU17} 6.3.7 (ii)).

\end{para}

\begin{para}\label{over1}
As above, let $u_1$ be the branch of $\log(1-x)$ and $T$ an indeterminate over $\cO_{\overline{X},1}$.
Then, by 1.1.1, we have an isomorphism $\cO_{\overline{X},y}^{\log}=\cO_{\overline{X},1}[u_1]\simeq \cO_{\overline{X},1}[T]$ of $\cO_{\overline{X},1}$-algebras under $(2\pi i)^{-1}u_1\leftrightarrow T$.
Consider an $\cO_{\overline{X},1}$-algebra homomorphism $\cO_{\overline{X},1}[T]\to\cO_{\overline{X},1}$, $T\mapsto 1-x$.

Let $\beta=(2\pi i)^{-1}\log(1-x)$.
Then, as $x\to 1$ in $\overline{X}$ along the real axis starting from $b$ over $0$ to $1$, $\exp(-\beta N_1)(\text{$F$ in \ref{naive1} (2)})$ converges to $F(0,0,-(2\pi i)^{-2}\zeta(2),\dots,-(2\pi i)^{-n}(c+\zeta(n)))$ in $\cD$ (\ref{hAlb3}), and hence the class of ($F$ in \ref{naive1} (2)) converges to the class $p_1$ (\ref{local1}) of the nilpotent orbit 
$$
(\sigma_1,\exp(\sigma_{1,\bC})F(0,0,-(2\pi i)^{-2}\zeta(2),\dots,-(2\pi i)^{-n}(c+\zeta(n))))
$$ 
in $\Gamma_{\sigma_1}\bs\cD_{\sigma_1}$.
We thus have an extension of the higher Albanese map over a simply connected neighborhood $\overline{X}_1$ of $1$ in $\overline{X}$ (Theorem \ref{t:FGS} (2)):
$$
\overline\varphi_1:\overline{X}_1\to \Gamma_{\sigma_1}\bs\cD_{\sigma_1}.
$$

This is a morphism in the category $\cB(\log)$.
The log structure on the source (resp.\ the target) is given by $1-x$ (resp.\ $q$).
The pullback of the universal log mixed Hodge structure on the target coincides with the log mixed Hodge structure on the source.

\end{para}

\begin{para}\label{log1}

By using log mixed Hodge theory, \ref{naive1} is described as follows.

Taking  the images of the nilpotent orbit in naive sense \ref{naive1} (1) and the \lq\lq higher Albanese map" \ref{naive1} (2), we have their real analytic extensions with boundary
$$
\overline\nu_1^{\log}, \;
\overline\varphi_1^{\log}:\overline{X}_1^{\log}\to(\Gamma_{\sigma_1}\bs\cD_{\sigma_1})^{\log}.
$$
Here, $\overline{X}_1^{\log}$ is similar to Example 1.1.1 over $x=1$, and
$(\Gamma_{\sigma_1}\bs\cD_{\sigma_1})^{\log}$ coincides with the moduli of nilpotent 
$i$-orbits $\Gamma_{\sigma_1}\bs\cD_{\sigma_1}^{\sharp}$ in the present situation (\cite{KNU} III Theorem 2.5.6).

Let $\tilde{\overline{X}}_1^{\log}$ be the universal covering of $\overline{X}_1^{\log}$.
The above maps are still lifted to
$$
\tilde{\overline\nu}_1^{\log}, \;
\tilde{\overline\varphi}_1^{\log}:\tilde{\overline{X}}_1^{\log}\to\cD_{\sigma_1}^{\sharp}.
$$
The boundary point $\tilde p_1$ in \ref{naive1} can be understood as the point $\tilde p_1=(z_1=0+i\infty)=(u_1=-\infty+i0)\in\tilde{\overline{X}}_1^{\log}$ (where $2\pi iz_1:=u_1$).
We have 
$(\exp(-(2\pi i)^{-1}\log(1-x)N_1)(\ref{naive1}\, (2)))(\tilde p_1)=F(0,0,-(2\pi i)^{-2}\zeta(2),\dots,-(2\pi i)^{-n}(c+\zeta(n)))$, and $\tilde{\overline\nu}_1^{\log}(\tilde p_1)=\tilde{\overline\varphi}_1^{\log}(\tilde p_1)\in\cD_{\sigma_1}^{\sharp}$ is the nilpotent $i$-orbit generated by $(N_1, F(0,0,-(2\pi i)^{-2}\zeta(2),\dots,-(2\pi i)^{-n}(c+\zeta(n))))$.

\end{para}

\begin{para}\label{eq'}

In order to describe the local structure near $x=\infty$, we take a local coordinate $\xi:=x^{-1}$.

By abuse of notation, let $F(\xi)$ be the pullback of the universal Hodge filtration by the composite $\varphi:X\to A_{X,\Gamma^{(n)}}\simeq\Gamma^{(n)}\bs\cD$ of higher Albanese map and the isomorphism in \ref{hAlb3}.

Since $d\log(x)=-d\log(\xi)$ and $-d\log(x-1)=d\log(\xi)-d\log(1-\xi)$, 
a connection equation in \ref{eq} now is
$$
dF(\xi)=\omega F(\xi),\quad
\omega:=(2\pi i)^{-1}\frac{d\xi}{\xi}(-N_0+N_1)+(2\pi i)^{-1}\frac{d\xi}{1-\xi}N_1.
$$
That is,
$$
da_{k-1,k}(\xi)=-(2\pi i)^{-1}\frac{d\xi}{\xi}\quad
(2\le k\le n),
$$
$$
da_{n,n+1}(\xi)=-(2\pi i)^{-1}\frac{d\xi}{\xi}-(2\pi i)^{-1}\frac{d\xi}{1-\xi},
$$
$$
da_{j,k}(\xi)=-(2\pi i)^{-1}a_{j+1,k}(\xi)\frac{d\xi}{\xi}\quad
(3\le k\le n+1,\;1\le j\le k-2).
$$

\end{para}

\begin{para}\label{sol'}

This system is solved by iterated integrals as before, and the solutions are
$$
a_{j,k}(\xi)=\frac{1}{(k-j)!}(-(2\pi i)^{-1}\log (\xi))^{k-j}\quad
(2\le k\le n,\;1\le j\le k-1),
$$
$$
a_{j,n+1}(\xi)=\frac{1}{(n+1-j)!}(-(2\pi i)^{-1}\log (\xi))^{n+1-j}
+(-(2\pi i)^{-1})^{n+1-j}l_{n+1-j}(\xi)\quad(1\le j\le n).
$$

\end{para}

\begin{para}\label{table'}

Table of solutions:
$$
\left(\CD
1\;\;&a_{1,2}\;\;&\ldots&a_{1,n}&a_{1,n+1}\\
0&1&\ddots&\vdots&\vdots\\
\vdots&0&\ddots\;\;&\;a_{n-1,n}\;\;&a_{n-1,n+1}\\
\vdots&\vdots&\ddots&1&a_{n,n+1}\\
0&0&\ldots&0&1
\endCD\right)
$$
$$
=\left(\CD
1\;\;&-(2\pi i)^{-1}\log(\xi)\;\;&\ldots\;\;&\frac{(-(2\pi i)^{-1}\log(\xi))^{n-1}}{(n-1)!}\;\;&
\frac{(-(2\pi i)^{-1}\log(\xi))^n}{n!}+(-(2\pi i)^{-1})^{n}l_n(\xi)\\
0&1&\ddots&\vdots&\vdots\\
\vdots&0&\ddots&-(2\pi i)^{-1}\log(\xi)&
\frac{(-(2\pi i)^{-1}\log(\xi))^2}{2!}+(-(2\pi i)^{-1})^{2}l_2(\xi)\\
\vdots&\vdots&\ddots&1&-(2\pi i)^{-1}\log(\xi)-(2\pi i)^{-1}l_{1}(\xi)\\
0&0&\ldots&0&1
\endCD\right).
$$

\end{para}

\begin{para}\label{local'}
Let now $\sigma_\infty=\bR_{\ge0}N_{\infty}$ with $N_{\infty}:=-N_0+N_1$ for $N_0$, $N_1$ in \ref{action}.
Let $F=F(-\alpha',\beta',\lambda_{2}',\dots,\lambda_{n}')$ be as in \ref{Hfilt}.
By Lemma \ref{l:MHSact} (i) (3), $(N_{\infty},F)$ satisfies the Griffiths transversality if and only if 
$\beta'=-\alpha',
\lambda_{2}'=\frac{(-\alpha')^2}{2!}, \dots,
\lambda_{n-1}'=\frac{(-\alpha')^{n-1}}{(n-1)!}$.
If this is the case, $(N_{\infty},F)$ generates a $\sigma_{\infty}$-nilpotent orbit, since admissibility and $\bR$-polarizability on $\gr^W$ trivially hold.
We describe the local structure of $\Gamma_{\sigma_{\infty}}\bs\cD_{\sigma_{\infty}}$ near the image $p_{\infty}$ of this nilpotent orbit.

Let $Y:=\{(q',\beta', \lambda_2',\dots,\lambda_n')\in \bC^{n+1}\,|\, \beta'=-\alpha',
\lambda_{2}'=\frac{(-\alpha')^2}{2!}, \dots,
\lambda_{n-1}'=\frac{(-\alpha')^{n-1}}{(n-1)!} \;\text{if}\; q'=0\}$ 
be the log manifold with the strong topology, with the structure sheaf of rings which is the inverse image of the sheaf of holomorphic functions on $\bC^{n+1}$, and with the log structure generated by $q'$. 
Then there is an open neighborhood $U$ of $(0,0,\dots,0,\lambda_n')$ in $\bC^{n+1}$ and an open immersion 
$$
Y\cap U \hookrightarrow \Gamma_{\sigma_{\infty}}\bs\cD_{\sigma_{\infty}}
$$ 
of log manifolds which sends $(q',\beta', \lambda_2',\dots,\lambda_n')\in Y\cap U$ with $q'\neq 0$ to the class of $F(-\alpha',\beta',\lambda_2',\dots,\lambda_n')$,  
where $\alpha'\in \C$ is such that $q'=e^{2\pi i\alpha'}$, and which sends $(0,0,\dots,0,\lambda_n')$ to $p_{\infty}$.

\end{para}

\begin{para}\label{naive'}

Near $x=\infty$, i.e., $\xi=0$, a nilpotent orbit in naive sense is
$$
(1)\hskip120pt
\exp((2\pi i)^{-1}\log(\xi)N_{\infty})F(0,0,\dots,0,\lambda_n^{\prime\,0})
\hskip120pt
$$
$$
=F\left(-(2\pi i)^{-1}\log(\xi),-(2\pi i)^{-1}\log(\xi),\frac{(-(2\pi i)^{-1}\log(\xi))^2}{2!},\dots,
\frac{(-(2\pi i)^{-1}\log(\xi))^n}{n!}+\lambda_n^{\prime\,0}\right),
$$
where $\lambda_n^{\prime\,0}=(-(2\pi i)^{-1})^nl_n(0)$.
The corresponding  \lq\lq higher Albanese map" (i.e., local version about $\infty$ of $\tilde\varphi$ in \ref{hAmap}) is
$$
(2)\hskip4pt
F\left(-(2\pi i)^{-1}\log(\xi),-(2\pi i)^{-1}\log(\xi)-(2\pi i)^{-1}l_1(\xi),\dots,\frac{(-(2\pi i)^{-1}\log(\xi))^n}{n!}+(-(2\pi i)^{-1})^nl_n(\xi)\right)\hskip4pt
$$
under the condition $l_j(0)=0$ ($1\le j\le n-1$).
These two are asymptotic when $\xi$ goes to the boundary point $b'$ described as follows.

Changing $\infty$ and $\xi$ into $0$ and $x$, respectively, $b'=(\infty,1)$ corresponds to the tangential boundary point $(0,1)$ of Deligne, i.e., $b'$ is the tangential base point over $\infty\in\overline{X}$ with tangent $v'\in T_\infty(\overline{X})=\Hom_\bC(m_\infty/m^2_\infty,\bC)$ defined by $v'(\xi)=1$.

This corresponds to our boundary point 
$b'=(y',a')$ with $y'=(\infty,h')\in\overline{X}^{\log}$ described as follows.
Let $u'$ be the branch of $\log(\xi)$ having real value on $\bR_{>0}$.
The argument function $h':M^{\gp}_{\overline{X},\infty}=\cO_{\overline{X},\infty}^{\times}\times \xi^\bZ \to {\bf S}^1$ is a group homomorphism sending $f\in\cO_{\overline{X},\infty}^{\times}$ to $f(\xi=0)/|f(\xi=0)|$ and $\xi$ to $v'(\xi)/|v'(\xi)|=1$, and the specialization $a':\cO^{\log}_{\overline{X},y'}=\bC\{\xi\}[u']\to \C$ is a ring homomorphism sending $\xi$ to $0$ and $u'$ to $a'(u')=\log(v'(\xi))=\log(1)=0$.

\end{para}

\begin{para}\label{over'}
As above, let $u'$ be the branch of $\log(\xi)$ and $T$ an  indeterminate over $\cO_{\overline{X},\infty}$.
Then, by 1.1.1, we have an isomorphism $\cO_{\overline{X},y'}^{\log}=\cO_{\overline{X},\infty}[u']\simeq \cO_{\overline{X},\infty}[T]$ of $\cO_{\overline{X},\infty}$-algebras under $(2\pi i)^{-1}u'\leftrightarrow T$.
Consider an $\cO_{\overline{X},\infty}$-algebra homomorphism $\cO_{\overline{X},\infty}[T]\to\cO_{\overline{X},\infty}$, $T\mapsto \xi$.

Let $\alpha'=(2\pi i)^{-1}\log(\xi)$.
Then, as $\xi\to0$, $\exp(-\alpha' N_\infty)(\text{$F$ in \ref{naive'} (2)})$ converges to $F(0,0,\dots,0,\lambda_n^{\prime\,0})$ in $\cD$ (\ref{hAlb3}), and hence the class of ($F$ in \ref{naive'} (2)) converges to the class $p_\infty$ (\ref{local'}) of the nilpotent orbit $(\sigma_\infty,\exp(\sigma_{\infty,\bC})F(0,0,\dots,0,\lambda_n^{\prime\,0}))$ in 
$\Gamma_{\sigma_\infty}\bs\cD_{\sigma_\infty}$.
We thus have an extension of the higher Albanese map over $\overline{X}_\infty$ (Theorem \ref{t:FGS} (2)):
$$
\overline\varphi_\infty:\overline{X}_\infty\to\Gamma_{\sigma_\infty}\bs\cD_{\sigma_\infty}.
$$

This is a morphism in the category $\cB(\log)$.
The log structure on the source (resp.\ the target) is given by $\xi$ (resp.\ $q$).
The pullback of the universal log mixed Hodge structure on the target coincides with the log mixed Hodge structure on the source.

\end{para}

\begin{para}\label{log'}

By using log mixed Hodge theory, \ref{naive'} is described as follows.

Taking the images of the nilpotent orbit in naive sense \ref{naive'} (1) and the \lq\lq higher Albanese map" \ref{naive'} (2), we have their real analytic extensions with boundary
$$
\overline\nu_\infty^{\log}, \;
\overline\varphi_\infty^{\log}:\overline{X}_\infty^{\log}\to(\Gamma_{\sigma_\infty}\bs\cD_{\sigma_\infty})^{\log}.
$$
Here, $\overline{X}_\infty^{\log}$ is like Example 1.1.1, and
$(\Gamma_{\sigma_\infty}\bs\cD_{\sigma_\infty})^{\log}$ coincides with the moduli of nilpotent 
$i$-orbits $\Gamma_{\sigma_\infty}\bs\cD_{\sigma_\infty}^{\sharp}$ in the present situation (\cite{KNU} III Theorem 2.5.6).

Let $\tilde{\overline{X}}_\infty^{\log}$ be the universal covering of $\overline{X}_\infty^{\log}$.
The above maps are still lifted to
$$
\tilde{\overline\nu}_\infty^{\log}, \;
\tilde{\overline\varphi}_\infty^{\log}:\tilde{\overline{X}}_\infty^{\log}\to\cD_{\sigma_\infty}^{\sharp}.
$$
The boundary point $b'$ in \ref{naive'} can be understood as the point $b'=(z'=0+i\infty)=(u'=-\infty+i0)\in\tilde{\overline{X}}_\infty^{\log}$ (where $2\pi iz':=u'$).
We have 
$(\exp(-(2\pi i)^{-1}\log(\xi)N_\infty)(\ref{naive'}\, (2)))(b')=F(0,0,\dots,0,\lambda_n^{\prime\,0})$, and 
$$
\tilde{\overline\nu}_\infty^{\log}(b')=\tilde{\overline\varphi}_\infty^{\log}(b')=
(\text{nilpotent $i$-orbit generated by $(N_\infty, F(0,0,\dots,0,\lambda_n^{\prime\,0}))$})\in\cD_{\sigma_\infty}^{\sharp}.
$$

\end{para}

\begin{para}\label{global}

For any $\sigma\in\Sigma$, $\Gamma_{\sigma}\bs\cD_\sigma\to\Gamma^{(n)}\bs\cD_\Sigma$ is a local homeomorphism.
This is analogously proved as \cite{KU09} Theorem A (iv).

Summing-up, we have a global extension over $\overline{X}$ of the higher Albanese map which is an isomorphism over its image:
$$
\overline\varphi:\overline{X}\overset\sim\to\overline\varphi(\overline{X})\subset  A_{X,\Gamma^{(n)},\Sigma}\simeq\Gamma^{(n)}\bs\cD_\Sigma.
$$

\end{para}

\begin{para}\label{prob}

To study analytic continuations and extensions of polylogariths in the spaces of nilpotent $i$-orbits $D_{\Sigma}^{\sharp}$, in the spaces of $\SL(2)$-orbits $D_{\SL(2)}$, and in spaces of Borel-Serre orbits $D_{\BS}$ is an interesting problem.
See \cite{KNU} for these extended period domains and their relations which are described as a fundamental diagram.

\end{para}

\setcounter{section}{0}
\def\thesection{\Alph{section}}

\section{Summary of a result of Deligne in \cite{D89}}\label{s:appen}
\medskip

We add here a summary of a result of Deligne in \cite{D89} for readers' convenience.

\begin{para}

Just as \ref{bsetting}--\ref{Gamma}, consider the situation $X:=\bP^1(\bC)\smallsetminus\{0,1,\infty\}\subset\overline{X}:=\bP^1(\bC)$.
Let $b:=(0,1)$ the \lq\lq tangential base point" over $0\in\overline{X}$ with tangent $1$.

Consider the quotient group $\Gamma$ of $\pi_1(X,b)$ as in \cite{D89} 16.14 (cf.\ \ref{Gamma}):
The inclusion $X\subset\bG_m(\bC)=\bC^{\times}$ induces $\pi_1(X,b)\to\pi_1(\bG_m(\bC),b)=\bZ(1)_B$ (suffix B means Betti, cf.\ \cite{D89}).
Let $K:=\Ker(\pi_1(X,b)\to\bZ(1)_B)$.
Let $\Gamma:=\pi_1(X,b)/[K,K]$ and $\Gamma_1:=K/[K,K]$.
Then, we have an exact sequence
$$1\to\Gamma_1\to\Gamma\to\bZ(1)_B\to1.$$

\end{para}

\begin{para} (\cite{D89} 16.15).
Let $\mu_0,\mu_1:\bZ(1)_B\to\Gamma$ be the monodromies around $0,1$, respectively.
Take a generator $u$ of $\bZ(1)_B$ (e.g.\ $u=2\pi i$), put $a_j=\mu_j(u)$ $(j=0,1)$.
Then, $\Gamma=\langle a_0,a_1\rangle$ with relation: conjugates of $a_1$ are commutative.

$\Gamma_1$ is a representation of $\bZ(1)_B$ with basis 
(conjugates of $a_1$) under the action
$\gamma \mapsto \mu_0(t)\gamma\mu_0(t)^{-1}$ ($\gamma\in\Gamma_1, t\in\bZ(1)_B$), i.e.,
$\Gamma_1=\bZ[\Z(1)_B]\cdot a_1$, where
$\sum_kc_k(a_0^ka_1a_0^{-k})=\sum_kc_k\cdot(2\pi i\cdot k)\cdot a_1$.

These are described as
$$
\Gamma_1=\bZ[\bZ(1)_B]\cdot a_1
\simeq\bZ[u,u^{-1}]\cdot\frac{du}{u},\quad
\Gamma=\bZ(1)_B\ltimes\Gamma_1,
$$
$$
\sum_kc_k(a_0^ka_1a_0^{-k})=\sum_kc_k\cdot(2\pi i\cdot k)\cdot a_1
\simeq\sum_kc_ku^k\frac{du}{u}
$$
(\cite{D89} 16.16).
Action of $\bZ(1)_B$ on $\Gamma_1$ is given by multiplication in $\bZ[\bZ(1)_B]=\bZ[u,u^{-1}]$.

\end{para}

\begin{para}

The descending central  series of $\Gamma$ induces a filtration on $\Gamma_1$:
$$
Z^N(\Gamma)\cap\Gamma_1=((u-1)^{N-1})\cdot\frac{du}{u}
\quad (N\geq 1)
$$
Let $\Gamma^{(N)}:=\Gamma/Z^{N+1}(\Gamma)$ and
$\Gamma_1^{(N)}:=\Image(\Gamma_1\to\Gamma^{(N)})$.
Then
$$
\Gamma_1^{(N)}=\bZ[u,u^{-1}]/(u-1)^N\cdot\frac{du}{u}.
$$

Put $u=e^v$ and hence $v=\log u$.
Then
$$
\bQ\otimes\Gamma_1^{(N)}=\bQ[u,u^{-1}]/(u-1)^N\cdot\frac{du}{u}
=\bQ[v]/(v^N)\cdot dv
$$
and we have
$$
\Gamma_1^{(N)}=
\left\{\sum_{k=0}^{N-1}c_k\exp(kv)dv\,\bigg|\,c_k\in\bZ\right\}.
$$

\end{para}

\begin{para}

As groups, identify
$$
\varphi:\bQ[v]/(v^N)\cdot dv\overset\sim\to\prod_1^N\bQ(n)_B:\quad
v^{n-1}dv=\frac{1}{n}dv^{n}\mapsto u^{\otimes n}.
$$
Then
$$
\sum_{k=0}^{N-1}c_k\exp(kv)dv\overset\varphi\mapsto
\hskip180pt
$$
$$
\hskip30pt
\sum_{k=0}^{N-1}c_k\left(\sum_{n=0}^{N-1}\frac{1}{n!}k^nu^{\otimes n}\right)\otimes u
=\sum_{n=1}^N\left(\sum_{k=0}^{N-1}c_k\frac{k^{n-1}}{(n-1)!}\right)u^{\otimes n}.
$$
Hence

\begin{sbprop} {\rm(\cite{D89} 16.17)}.\quad
$(n-1)!\cdot{\rm pr}_n\circ\varphi(\Gamma_1^{(N)})=\bZ(n)_B$.
\end{sbprop}

\end{para}

\begin{para} (\cite{D89} 16.12). 
Define a Lie algebra action of $\bQ(1)$ on $\prod_1^N\bQ(n)$ by 
$$
a\ast(b_1,b_2,\dots,b_N)=(0,ab_1,\dots,ab_{N-1}),
$$
and $\bQ(1)\ltimes\prod_1^N\bQ(n)$ the associated semi-direct product of Lie algebra.

Let $\mu_0, \mu_1:\bQ(1)\to\bQ(1)\ltimes\prod_1^N\bQ(n)$ be morphisms of Lie algebras such that $\mu_0$ is the identity onto the first factor $\bQ(1)$ and $\mu_1$ is the identity onto  the factor $\bQ(1)$ in the product $\prod_1^N\bQ(n)$.

By abuse of notation, let $\mu_0,\mu_1:\bQ(1)\to \bQ\otimes\Lie \Gamma^{(N)}$.
Then there exists a unique Lie algebra isomorphism respecting each $\mu_0$, $\mu_1$:
$$
\bQ(1)\ltimes\prod_1^N\bQ(n) \overset\sim\to\bQ\otimes\Lie \Gamma^{(N)}
=\bQ(1)\ltimes(\bQ\otimes\Lie\Gamma_1^{(N)})
$$
which is given by $\mu_0$ and
$\nu_n:=(\operatorname{ad}\mu_0)^{n-1}(\mu_1)$ $(1\le n\le N)$.

\end{para}

\begin{para}\label{itTate}

Let $\operatorname{Lie} U_{\text{DR}}^{(N)}$ be the de Rham realization of iterated Tate motive in \cite{D89} 16.13.
Let $e_\alpha:=\mu_\alpha(1)\in\operatorname{Lie} U_{\text{DR}}^{(N)}$ 
($1=\exp(2\pi i)\in\bQ(1)_{\text{DR}}$, $\alpha=0,1$).

Take  coordinates $(u,(v_n)_{1\le n\le N})$ of $U_{\text{DR}}^{(N)}$ as follows:
$$
(u,(v_n)_n)\mapsto
\exp(ue_0)\exp\left(\sum_{n=1}^{N} v_n(\Ad e_0)^{n-1}(e_1)\right).
$$ 

\begin{sblem}  {\rm(\cite{D89} 19.3.1).}
Let $z\in\bC^{\times}\smallsetminus\bR_{\ge1}$.
The end point of the image in $U_{\text{DR}}^{(N)}(\bC)$ of the line segment from $(0,z)$ to $z$ has coordinates $u=0$, $v_n=-l_n(z)$.
\end{sblem}

\noindent
{\it Proof.}
Let $z_1,z_2\in\bC^{\times}\smallsetminus\bR_{\ge1}$.
Take a path from $z_1$ to $z_2$, and take an iterated integral $I_{z_1}^{z_2}$ of 
$$
dI(t)=\left(\frac{dt}{t}e_0+\frac{dt}{t-1}e_1\right)\cdot I(t)
$$
for $I(t)=1+ue_{0}+\sum_{n}v_{n}(\Ad e_{0})^{n-1}(e_{1})$.
Note
$e_0\ast e_0=e_0$, $e_0\ast(\Ad e_0)^{n-1}(e_1)=(\Ad e_0)^n(e_1)$\;\;$(1\le n\le N)$,
$e_1\ast e_0=0$, $e_1\ast e_1=e_1$, $e_1\ast(\Ad e_0)^{n-1}(e_1)=0$\;\;$(2\le n\le N)$.

The corresponding differential equation is 
$$
du=\frac{dt}{t}, \quad 
dv_1=\frac{dt}{t-1}, \quad
dv_n=v_{n-1}\frac{dt}{t}.
$$
Take $I(z_1)=\text{identity}\in U_{\text{DR}}^{(N)}(\bC)$ as an initial condition and consider $z_2$ as a variable.

If $z_1$ is a tangential base point $(0,\tau)$ (\cite{D89} Section 15), replace the initial condition by
$$
I(t)\exp\left(-\log\left(\frac{t}{\tau}\right)\right)\to\text{identity}\quad
\text{as $t\to0$}.
$$
For the line segment from $(0,z)$ to $z$, we have
$$
u=\log\left(\frac{t}{z}\right), \quad
v_n=-l_n(t).
\qed
$$

\end{para}
\bigskip

{\bf Acknowlegements.}
The author thanks Kazuya Kato and Chikara Nakayama for joint research.
This note grew up in the preparation for and during the workshop \lq\lq Hodge theory and algebraic geometry" held at the end of August, 2017, Tokyo Denki University.
The author thanks Tetsushi Ito, Satoshi Minabe, Taro Fujisawa, and Atsushi Ikeda for the good occasion and helpful discussions. 
The author thanks the referee for careful reading and valuable comments.

S.\ Usui was 
partially supported by JSPS Grants-in-Aid for Scientific Research (C) 17K05200.

\bigskip

\noindent
{\rm Sampei USUI
\\
Graduate School of Science
\\
Osaka University
\\
Toyonaka, Osaka, 560-0043, Japan}
\\
{\tt usui@math.sci.osaka-u.ac.jp}

\end{document}